\documentclass[journal,twoside,web]{ieeecolor}
\usepackage{generic}
\usepackage{cite}
\usepackage{amsmath,amssymb,amsfonts}
\usepackage{algorithmic}
\usepackage{graphicx}
\usepackage{algorithm,algorithmic}
\usepackage{hyperref}
\usepackage{textcomp}
\usepackage{subfig}
\usepackage{epstopdf}
\usepackage{threeparttable}
\newtheorem{theorem}{Theorem}
\newtheorem{corollary}{Corollary}
\newtheorem{definition}{Definition}
\newtheorem{lemma}{Lemma}
\newtheorem{proposition}{Proposition}
\newtheorem{assumption}{Assumption}

\newtheorem{remark}{Remark}

\usepackage{xcolor}
\usepackage{colortbl}
\usepackage{tikz}
\newcommand{\ALGNAME}{AHEAD}

\hypersetup{
  colorlinks=false,
  pdfborder={0 0 0}
}

\definecolor{subcolor}{RGB}{0,0.541,0.855}

\def \colorbf #1{\textcolor{subcolor}{#1}}

\def \reNYC   #1{\textcolor{blue}{#1}}

\usepackage{ifthen}
\def\ishidecomments{no}  
\def\authortextcolor{red} 
\newcommand{\jm}[2][]{%
   \ifthenelse{\equal{#1}{}}
   {\textcolor{\authortextcolor}{#2}~\ignorespaces}
      {\ifthenelse{\equal{\ishidecomments}{yes}}
         {\textcolor{\authortextcolor}{#2}~\ignorespaces}
         {\textcolor{\authortextcolor}{#2}~\textcolor{cyan}{[\textbf{jm}:~#1]}~\ignorespaces}
      }
}

\def\BibTeX{{\rm B\kern-.05em{\sc i\kern-.025em b}\kern-.08em
    T\kern-.1667em\lower.7ex\hbox{E}\kern-.125emX}}

\begin{document}

\title{Hessian-Free Distributed Bilevel Optimization via Penalization with Time-Scale Separation}



\author{Youcheng Niu$^\dagger$, Jinming Xu$^{\dagger,\star}$, Ying Sun$^\ddagger$,
Li Chai$^\dagger$, and Jiming Chen$^\dagger$
\thanks{This work was supported in part by National Natural Science Foundation of China (NSFC)
under Grant 62373323, Grant W2511069, Grant 62088101, and Grant U244124, and in
part by Zhejiang Province Natural Science Foundation of China under Grant LZ24F030006. }
\thanks{$^\dagger$Y. Niu, J. Xu, L. Chai, and J. Chen are with the State Key Laboratory of Industrial Control Technology and the College of Control Science and Engineering, Zhejiang University, Hangzhou 310027, China.}
\thanks{$^\ddagger$Y. Sun is with the School of Electrical Engineering and Computer Science, Pennsylvania State University, PA 16802, USA.}
\thanks{$^\star$Correspondence to Jinming Xu ({\tt jimmyxu@zju.edu.cn})}
}

\maketitle

\begin{abstract}
This paper considers a class of  distributed bilevel optimization (DBO) problems with a coupled inner-level subproblem. Existing approaches typically rely on hypergradient estimations involving computationally expensive Hessian evaluation.
To address this, we approximate the DBO problem as a minimax problem by properly designing a penalty term that enforces both the constraint imposed by the inner-level subproblem and the consensus among the decision variables of agents.
Moreover, we propose a loopless distributed algorithm, \ALGNAME, that employs multiple-timescale updates to solve the approximate problem asymptotically without requiring Hessian computation.
Theoretically, we establish  sharp convergence rates for  nonconvex-strongly-convex settings and for distributed minimax problems as special cases.
Our analysis reveals a clear dependence of convergence performance on node heterogeneity, penalty parameters, and network connectivity, with a weaker assumption on heterogeneity that only requires bounded gradients at the optimum. Numerical experiments corroborate our theoretical results.

\end{abstract}


\section{Introduction}
\label{sec:introduction}
\IEEEPARstart{D}{istributed} bilevel optimization (DBO) has attracted growing attention in recent years, driven by the increasing demand for large-scale, decentralized computation and its broad range of potential applications  such as Stackelberg games \cite{fabiani2021local}, meta-learning \cite{MAML}, hyperparameter optimization \cite{bertrand2020implicit}, distributed minimax problems \cite{razaviyayn2020nonconvex,notarnicola2018duality}, multiple-agent coordination \cite{carnevale2024nonconvex},  and  robot control \cite{wang2024imperative}.
In distributed bilevel optimization, multiple nodes collaborate to minimize a sum of local outer-level objectives, where these objectives depend on parameters typically determined by feasible solutions to corresponding inner-level subproblems. These local objectives and computations are distributed across a network of nodes. The problem can be formulated as:
$
\min_{x \in \mathbb{R}^n} \Phi(x) = \frac{1}{m}\sum_{i=1}^m f_i(x, y^*(x)) \; \text{s.t.} \; y^*(x) \in \mathcal{S}(x),
$
where $f_i: \mathbb{R}^n \times \mathbb{R}^r \rightarrow \mathbb{R}$ is the outer-level objective  known only to node $i$ and \( \mathcal{S}(x) = \arg\min_{y \in \mathbb{R}^r} g(x,y) \) represents the set of optimal solutions to the inner-level subproblem associated with the objective \( g(x,y) \).
This distributed bilevel framework naturally models interactions between tasks and nodes in practical applications. For example, in energy management, utility companies set electricity prices \( x \) based on consumer consumption patterns \( y^*(x)\), while consumers adjust their consumption in response \cite{maharjan2013dependable}.

Distributed optimization offers several key advantages, such as scalability, privacy preservation, and efficient use of computing resources \cite{nedic2018distributed, wang2022gradient}. Previous works have primarily focused on the single-level problems in the form of
$
\min_{x \in \mathbb{R}^n} \frac{1}{m}\sum_{i=1}^{m} f_i(x),
$
where \( f_i(x) \) represents the local objective at node \(i\). Classical methods, including distributed gradient descent (DGD) \cite{nedic2017achieving, li2018distributed,lian2017can}, gradient tracking (GT) \cite{xu2015augmented, sun2022distributed}, and primal-dual methods \cite{shi2014linear, alghunaim2020linear}, have established the foundation for solving such single-level optimization problems in distributed systems. Initially focused on convex objectives, these methods have more recently been extended to handle nonconvex objectives. Since many practical applications involve constraints, researchers have turned to distributed constrained optimization (DCO), developing various methods to address local or coupled inequality and equality constraints \cite{li2020distributed, chang2014distributed, towfic2014adaptive}. These approaches primarily consist of penalty-based and Lagrangian-based methods. For instance, Towfic \textit{et al.} \cite{towfic2014adaptive} investigate distributed penalty-based methods for both local inequality and equality constraints, which are favored for their simplicity and ease of implementation. In contrast, Li \textit{et al.} introduce dual variables within Lagrangian-based methods to explicitly manage global inequality constraints \cite{li2020distributed}. Notably, the DBO  problem can reduce to the abovementioned DCO problem when the feasible solution set \(\mathcal{S}(x)\) can be explicitly determined. However, this is often not  the case, as the complexity of the inner-level subproblem in many practical applications renders the singe-level distributed methods invalid.

Recent research has made significant strides in addressing the DBO problem, proposing efficient algorithms based on approximate implicit differentiation (AID) methods that improve performance and scalability, particularly when the feasible set \( \mathcal{S}(x) \) is a singleton and does not have an explicit expression \cite{yang2022decentralized, chen2023decentralized, niu2023loopless, dong2023single, niu2023distributed, kong2024decentralized,hong2023two}. For example, the work in \cite{chen2023decentralized} tackles a stochastic setting by proposing a double-loop  algorithm based on the DGD method, incorporating additional computation loops to solve the inner-level subproblem and to approximate the Hessian-inverse-vector product for hypergradient estimation. However, this results in high computational costs for estimating the hypergradient. To mitigate this complexity, \cite{niu2023loopless} proposes a loopless distributed algorithm that leverages one-step approximation techniques to address problems with personalized inner-level subproblems. After that, several other loopless distributed algorithms have been proposed recently~\cite{niu2023distributed, kong2024decentralized, dong2023single, zhu2024sparkle}. Note that early works such as~\cite{niu2023loopless, chen2023decentralized, yang2022decentralized} rely on a strong bounded gradient assumption, and do not clarify how the convergence is  affected by node heterogeneity, unlike the treatment in single-level optimization \cite{lian2017can}, i.e., the variance  in the objective functions across  nodes. To address this limitation, \cite{niu2023distributed} introduces a unified loopless framework that combines the DGD and GT methods and provides a theoretical analysis of both outer- and inner-level heterogeneity, elucidating their effects on convergence in personalized settings. For DBO problems with coupled inner-level subproblems, \cite{kong2024decentralized} analyzes the transient iteration complexity using DGD methods, while \cite{dong2023single} designs a loopless GT-based algorithm specifically for deterministic scenarios. Importantly, all these  methods~\cite{yang2022decentralized, chen2023decentralized, niu2023loopless, dong2023single, niu2023distributed, kong2024decentralized, zhu2024sparkle} estimate the hypergradient, i.e., \( \nabla \Phi(x) \), which entails costly computation due to the involvement of Hessian-related information.
 Moreover, existing studies on DBO \cite{niu2023distributed, kong2024decentralized} generally impose stronger second-order conditions for heterogeneity analysis, unlike single-level optimization \cite{lian2017can}. It remains unclear whether second-order heterogeneity is truly essential for heterogeneity analysis in DBO.

Motivated by the interpretation of \( y^*(x) \in \mathcal{S}(x) \) as a potential constraint, pioneering works \cite{liu2022bome, kwon2023fully, yao2024constrained} employ value functions (the optimum value of the inner-level objective, given the outer-level variable) to reformulate bilevel optimization as equivalent constrained problems (known as  the value-function-based problem (VP)) by constraining the inner-level function's value to be less than or equal to the value function~\cite{ye2023difference}. By leveraging the value function to implicitly capture the derivative of the inner-level solution, these works bypass the need of Hessian estimation, achieving a convergence  rate of \( \mathcal{O}\left(\frac{1}{K^{2/3}}\right) \) for nonconvex-strongly-convex cases. However, these methods \cite{liu2022bome, kwon2023fully, ye2023difference, yao2024constrained} are not directly applicable in distributed settings. In particular, in the distributed context, local value functions  depend on a global inner-level solution, making their  gradient estimation involve the implicit derivative of the global solution and Hessian estimation, which differs from those in centralized scenarios. Moreover, the absence of direct access to global information in sparse networks leads to consensus errors, degrading constraint satisfaction in the VP problem and potentially hindering convergence. These factors make solving the DBO problem in a Hessian-free manner with guaranteed convergence particularly challenging.
\textbf{Summary of Contributions.} To tackle  the above challenges, this work proposes a new computationally efficient distributed algorithm for  the DBO problem, while  providing a sublinear
convergence guarantee in  nonconvex–strongly-convex settings, along with a heterogeneity analysis. The main contributions  are summarized as follows:
\begin{itemize}
   \item  We propose a new penalty-based distributed optimization algorithm, termed \ALGNAME, which is Hessian-free and features a loopless structure, setting it apart from previous approaches \cite{yang2022decentralized, chen2023decentralized, niu2023loopless, dong2023single, niu2023distributed, kong2024decentralized, zhu2024sparkle}.
     In particular, by introducing proper penalty terms, we propose a  single-level approximation with  minimax formulations for the DBO problem, without the need to estimate the implicit derivative of the inner-level solution  in a distributed setting.
   Different from \cite{liu2022bome, kwon2023fully, ye2023difference, yao2024constrained, li2020distributed, chang2014distributed, towfic2014adaptive}, the penalty term    enforces both the constraint
induced by the inner-level subproblem as well as the consensus among the decision variables. Moreover, we employ multiple-timescale
   updates according to the subproblem's  properties to ensure algorithm's convergence, while maintaining a loopless structure.

    \item Theoretically, we explicitly characterize the dynamics related to the  errors induced by the penalization and node heterogeneity, as well as the impact of the penalty parameter on convergence performance. In particular, we establish a  rate of $\mathcal{O}\big( \frac{\kappa ^4}{K^{1/3}}\!\!+\!\!\frac{1}{(1-\rho)^2}(\frac{{\kappa ^2b_{f}^{2}}}{K^{2/3}}\!+\!\frac{\kappa ^4b_{g}^{2}}{K^{1/3}}) \big)$  for nonconvex-strongly-convex cases, which clearly shows how the out- and inner-level  heterogeneity  $b_f^2$,  $b_g^2$, network connectivity \(\rho\), and  the condition number  \(\kappa\) effect the convergence. Besides, we provide a tight heterogeneity analysis by exploring the relationship between the estimates for the inner-level solution and their penalized counterparts, which relies solely on the first-order  heterogeneity for both level and  its boundedness at the optimum, distinguishing it from existing works on DBO \cite{niu2023distributed, kong2024decentralized}.
    Moreover, we extend our result to distributed minimax problems as a special case, achieving a faster rate of  $\mathcal{O}\big(\! \frac{\kappa^2}{ K^{\frac{2}{3}}  }\!\!+\! \!\frac{\kappa^2 b_{f}^{2}}{\left( 1-\rho \right) ^2  K^{\frac{2}{3}} } \big)$.
\end{itemize}
We note that a recent parallel work \cite{wang2024fully} develops a Hessian-free distributed DSGDA-GT algorithm based on the gradient tracking scheme to solve DBO in stochastic settings. This work achieves a convergence rate of \(\mathcal{O}(\frac{1}{K^{1/7}})\) given that the inner-level variables are initialized close to optimal solutions, that higher-order smoothness holds for the outer-level functions and that the gradient of the inner-level functions, $\|\nabla_xg_i(x,y)\|$, is bounded.
In contrast, we approach the problem from a penalty-based perspective, with an emphasis on heterogeneity analysis. We do not require additional initial conditions and provide a tight convergence result under weaker assumptions, revealing the clear dependence of the convergence performance on several factors such as node heterogeneity,  network connectivity, and the condition number.
\\
\noindent \textbf{Notation.} Throughout this paper, vectors are assumed to be column vectors unless otherwise stated. \( I_n \) denotes the \( n \times n \) identity matrix and \( 1_m \) the \( m \times 1 \) all-ones vector. Given a vector \( x \), \( \|x\| \) denotes its Euclidean norm. The vector formed by stacking column vectors \( x_1, \cdots, x_m \) on top of each other is denoted as \( {\rm col}\{x_i\}_{i=1}^m \). For simplicity, the notation $\{x_i\}$ without the script indicates \( x_1, \cdots, x_m \).    For a graph \( \mathcal{G} = \{\mathcal{V}, \mathcal{E}\} \), we use the symbol \( \mathcal{N}_i \) to denote the neighbors of node \( i \), with the node \( i \) itself being included, i.e., \( i \in \mathcal{N}_i \). Given a weight matrix $W \in \mathbb{R}^{n \times n}$ of a graph,  we denote its Kronecker product with the identity matrix $I_n$ by $\mathcal{W}_n$, i.e., $\mathcal{W}_n=W \otimes {I_n}$. For a function \( f(x, y): \mathbb{R}^n \times \mathbb{R}^r \to \mathbb{R} \), \( \nabla_x f \) and \( \nabla_y f \) represent the partial gradients w.r.t. \( x \) and \( y \), respectively, while \( \nabla_{xy}^2 f \) and \( \nabla_{yy}^2 f \) denote the Jacobian and Hessian matrices, respectively.

\section{Problem Formulation}
Consider an undirected and connected  network $\mathcal{G}=\{\mathcal{V},\mathcal{E}\}$ with  $\mathcal{V}=\{1,\cdots, m\}$ being the node set and $\mathcal{E}
\subset \mathcal{V} \times \mathcal{V} $ the edge set,  where $m$ nodes cooperatively solve  a class of DBO problems with a coupled inner-level subproblem as follows:
\begin{equation}\label{EQ-P1}
\begin{aligned}
&\underset{x\in \mathbb{R} ^n}{\min}\Phi \left( x \right) \triangleq\frac{1}{m}{\sum_{i=1}^m{f_i\left( x,y^*\left( x \right) \right)}}, \;  \\
&\mathrm{s}.\mathrm{t}. \; y^*\left( x \right) =\mathrm{arg}\underset{y\in \mathbb{R} ^r}{\min}g\left( x,y \right) \triangleq\frac{1}{m}\sum_{i=1}^m{g_i\left( x,y \right)},
\end{aligned}
\end{equation}
where $f_i:\mathbb{R} ^n\times \mathbb{R} ^r\rightarrow \mathbb{R}$ is the  outer-level function  and  $g_i:\mathbb{R} ^n\times \mathbb{R} ^r\rightarrow \mathbb{R}$ represents the inner-level function.  In the problem \eqref{EQ-P1}, both \( f_i \) and \( g_i \) are local to node \( i \) and cannot be accessed by other nodes.
For convenience, we denote $f\left( x,y^*\left( x \right) \right)=(1/m){\sum_{i=1}^m{f_i\left( x,y^*\left( x \right) \right)}}$. Now, we make the following assumptions.
\begin{assumption} [\textbf{Outer-level functions}] \label{ASS-outer-level}
We assume that
i) $f_i(x,y)$ is continuously differentiable and $L_{f,1}$-smooth jointly in $(x,y)$;
ii) $
\left\| \nabla _y f_i \left( x,y^*(x) \right) \right\|
$ is bounded by $C_{f,y}$.
\end{assumption}
\begin{assumption}[\textbf{Inner-level functions}]\label{ASS-inner-level}
We assume that
 i) $g_i(x,y)$ is $\mu_g$-strongly convex in $y$; ii)  $g_i(x,y)$ is continuously differentiable and $L_{g,1}$-smooth jointly in $(x,y)$; iii)  $\nabla ^2 g_i$ is $L_{g,2}$-Lipschitz continuous in $(x,y)$.

\end{assumption}
\begin{assumption}[\textbf{Network connectivity}]
\label{ASS-network}
The  network \(\mathcal{G}\) is undirected and connected. The weight matrix $W$ of  \(\mathcal{G}\) is doubly stochastic such that: i) \(w_{ij} = w_{ji} > 0\) if and only if \((i, j) \in \mathcal{E}\); otherwise, \(w_{ij} = 0\); ii) $\sum_{j \in \mathcal{N}_i} w_{ij}=1$, $
\forall i \in \mathcal{V}$.
\end{assumption}

Assumptions \ref{ASS-outer-level} and \ref{ASS-inner-level} characterize the smoothness of the outer- and inner-level functions, which are commonly adopted in bilevel optimization literature \cite{ghadimi2018approximation,kong2024decentralized,chen2023decentralized,yang2022decentralized}. Smoothness assumptions are pervasive and indispensable in bilevel optimization literature, as they lay the theoretical groundwork for analyzing the convergence behavior of optimization algorithms, deriving error bounds, and ensuring the well-posedness of the optimization problem.   The strong convexity of the inner-level function \( g_i \) in \( y \)  guarantees the uniqueness of \( y^*(x) \) for any \( x \) and ensures the differentiability of the function \( \Phi(x) \). Notably, Assumption \ref{ASS-outer-level}(ii) requires only the boundedness of \(\|\nabla_y f_i\|\) at the point \(y^*(x)\), which is a weaker condition compared to those in \cite{kwon2023fully,yang2022decentralized}, where the boundedness is required for all \(y\). Moreover, Assumption \ref{ASS-network} is standard for ensuring network connectivity in distributed optimization \cite{xu2015augmented, alghunaim2020linear}. Note that \(\rho \triangleq \| W - \frac{1_m 1_m^{\rm{T}}}{m} \|^2 \in [0, 1)\) holds under Assumption \ref{ASS-network} \cite{xu2015augmented}.

\subsection{Penalty-based Single-Level Approximation}
To address the DBO problem in a computationally efficient manner, we seek to avoid the  direct estimation of the optimal solution to the inner-level subproblem and the computation of  the implicit derivative of the inner-level solution. To achieve this, we introduce an auxiliary variable \(z\) as a proxy for \(y^*(x)\), and allow each node \(i\) to maintain local copies of \(x_i\), \(y_i\), and \(z_i\). We then propose a  single-level approximate problem  by  leveraging  penalty-based techniques, as  follows:
\begin{align}
&\min_{\left\{ x_i \right\} ,\left\{ y_i \right\}} \max_{\left\{ z_i \right\}} \frac{1}{m}\sum_{i=1}^m{\left( f_i(x_i,y_i)+\lambda (g_i(x_i,y_i)-g_i(x_i,z_i)) \right)}
\nonumber\\
&\;+\!\frac{1}{m}\!\sum_{i=1}^m\!{\sum_{j\in \mathcal{N} _i\backslash\{i\}}\!\!\!\!{\big( \frac{1}{2\alpha}\!\left\| \sqrt{w_{ij}}\!\left( x_i\!-\!x_j \right) \right\| ^2\!\!+\!\frac{1}{2\beta}\!\left\| \sqrt{w_{ij}}\!\left( y_i\!-\!y_j \right) \right\| ^2 \! \big)}}
\nonumber\\
&\;-\!\frac{1}{m}\!\sum_{i=1}^m\!{\sum_{j \in \mathcal{N} _i\reNYC{\backslash}\{i\}}\!\!{\big( \frac{\lambda}{2\gamma}\left\| \sqrt{w_{ij}}\!\left( z_i-z_j \right) \right\| ^2 \big)}},  \label{EQ-P-penalty-adptive}
\end{align}
where $\lambda>0$ is the  penalty parameter for the constraint imposed by the inner-level subproblem, while the coefficients  \(\frac{1}{\alpha}\), \(\frac{1}{\beta}\) and \(\frac{\lambda}{\gamma}\) are the  penalty amplitudes for the  consensus constraints.  Here, \(\alpha>0\), \(\beta>0\) and \(\gamma>0\) are  multiple-timescale parameters to be determined later. In  \eqref{EQ-P-penalty-adptive},  both the constraint  \( y \in y^*(x) \)  induced by the inner-level subproblem and the consensus constraint are jointly penalized leveraging the problem's inherent properties. In particular, the maximization subproblem ensures that the auxiliary variables \( \{z_i\} \) can effectively  approximate the inner-level solution \( y^*(x) \). In the minimization subproblem over $\{x_i\}$ and $\{y_i\}$, the penalty term \( \lambda (g_i(x_i, y_i) - g_i(x_i, z_i)) \), combined with the consensus penalty, properly enforces local consistency between \(\{ y_i\} \) and \( \{z_i\}\). This allows \( \{y_i\} \) to remain close to the inner-level solution while  minimizing the outer-level objectives, with the penalty parameter $\lambda$ controlling the penalty gap at both the outer and inner levels.
Notably, we introduce an additional \(\lambda\)-dependent consensus penalty term associated with the auxiliary variables  \( \{z_i\} \) in the objective function.
This term enhances the consistency  among the local estimates $\{z_i\}$ across the network, thereby improving their accuracy in approximating the optimal inner-level solution $y^*(x)$. Such enhancement is also essential to ensure that  the approximation quality of $\{z_i\}$ remains consistent with that of $\{y_i\}$, as required by the  \(\lambda\)-dependent penalty term $\lambda (g_i(x_i, y_i) - g_i(x_i, z_i))$ in the distributed setting.

\begin{remark}[\textbf{Hessian-free property}] It should be noted that previous DBO methods \cite{yang2022decentralized, chen2023decentralized, niu2023loopless, dong2023single, niu2023distributed, kong2024decentralized, zhu2024sparkle} typically need to directly estimate the hypergradient $\nabla \Phi(x)$ as follows\cite{ghadimi2018approximation}:
\begin{equation}
\begin{aligned} \label{EQ-hypergradient}
\nabla \Phi \left( x \right) =\nabla _xf\left( x,y^*\left( x \right) \right)+\nabla y^*(x)\nabla _yf\left( x,y^*\left( x \right) \right),
\end{aligned}
\end{equation}
which relies on the implicit derivative \(\nabla y^*(x) = -\nabla _{xy}^{2}g\left( x, y^*(x) \right) \left[ \nabla _{yy}^{2}g\left( x, y^*(x) \right) \right]^{-1}\) that involves the evaluation of Hessians. In contrast, we introduce the  penalty term \( \lambda (g_i(x_i, y_i) - g_i(x_i, z_i)) \)  to enforce the constraint imposed by the inner-level subproblem, and thus the  single-level approximate problem  \eqref{EQ-P-penalty-adptive} allows us to bypass  the direct estimation for \(\nabla y^*(x)\) and  Hessian  in the distributed setting.

\end{remark}

\section{Algorithm  and Main Results }
In this section, we propose a loopless distributed algorithm, termed {\ALGNAME}, in a Hessian-free manner based on the  approximation in \eqref{EQ-P-penalty-adptive}, and present the main  results.

Before we present the algorithm and main results, we first introduce the
following notations for brevity:
\begin{equation}
\begin{aligned}
&x^k \triangleq\mathrm{col}\{x_i^k\}_{i=1}^m,\; y^k \triangleq  \mathrm{col}\{y_i^k\}_{i=1}^m, \;z^k \triangleq  \mathrm{col}\{z_i^k\}_{i=1}^m, \\
&\bar{x}^k \triangleq\frac{1}{m}\sum_{i=1}^m x_i^k, \;\bar{y}^k  \triangleq  \frac{1}{m}\sum_{i=1}^m y_i^k, \;\bar{z}^k \triangleq  \frac{1}{m}\sum_{i=1}^m z_i^k, \\
\end{aligned}
\end{equation}
where  \(x_i^k\), \(y_i^k\) and \(z_i^k\)  denote the estimates of the variables \(x\), \(y\) and \(z\) at iteration \(k\) for node \(i\), respectively; and we denote an  auxiliary penalty function by $p\left( x,y,z;\lambda \right) \triangleq\frac{1}{m}\sum_{i=1}^m{\left( f_i\left( x,y \right) +\lambda q_i\left( x,y,z \right) \right)}$ with $q_i\left( x,y,z \right) \triangleq g_i(x,y)-g_i(x,z)$ as well as  its corresponding solution with respect to $y$ by $y^*\left( x;\lambda \right)$ (referred to as the  penalized inner-level solution), i.e.,  $y^*\left( x;\lambda \right) \in \arg\min _y\max _zp\left( x,y,z;\lambda \right)$, where $y^*\left( x;\lambda \right)$ exists and is unique under an appropriate choice of $\lambda$ (c.f., Lemma \ref{LE-inner-penalty-error}).

We also use the following compact notations for brevity:
\begin{equation}
\begin{aligned}
\nabla _x G(x^k,y^k)&\triangleq{\rm{col}}\{\nabla _x g_i(x_i^k,y_i^k)\}_{i=1}^m,\\
 \nabla _x G(x^k,z^k)&\triangleq{\rm{col}}\{\nabla _x g_i(x_i^k,z_i^k)\}_{i=1}^m,\\
\nabla _x Q (x^k,y^k,z^k)&\triangleq {\rm{col}}\{\nabla _x q_i(x_i^k,y_i^k, z_i^k)\}_{i=1}^m, \nonumber \\
\end{aligned}
\end{equation}
with  $\nabla _x F(x^k,y^k)$, $\nabla _y F(x^k,y^k)$, and $ \nabla _y G(x^k,y^k)$ defined analogously.  Also, we define the following constants induced by the Lipschitz continuity and strong convexity (c.f., Assumptions \ref{ASS-outer-level} and \ref{ASS-inner-level}):
\begin{align}
U_{\lambda}^2 &\triangleq L_{f,1}^{2}\!+\lambda ^2L_{g,1}^{2}, \;
L_{\lambda}\triangleq L_{f,1}+\lambda L_{g,1}, \nonumber\\
\mu _{\lambda}&\triangleq\frac{\lambda \mu_g}{2},\;
L_{y^*}
\triangleq\frac{L_{g,1}}{\mu _g} , \;L_{y^*,\lambda}\triangleq
\frac{2L_{\lambda}}{\lambda \mu _g}, \label{EQ-constant-symbols}
\\
L&\triangleq\! ( L_{f,1}\!+\!\frac{(L_{f,1}\!+\!C_{f,y})L_{g,2}}{\mu _g}\!+\!\frac{C_{f,y}L_{g,2}^{2}}{\mu _{g}^{2}}\! ) ( 1\!+\!L_{y^*} \!), \nonumber
\\
w_{\gamma}&\triangleq\frac{\mu _gL_{g,1}}{2(\mu _g+L_{g,1})}, \; w_{\beta}\triangleq\frac{\mu _{\lambda}L_{\lambda}}{2(\mu _{\lambda}+L_{\lambda})}, \nonumber
\end{align}
and the following constants associated with the errors induced by the penalty term in \eqref{EQ-P-penalty-adptive}:
\begin{align}
C_{\mathrm{in}}\triangleq\!\frac{2C_{f,y}}{\mu _g}, C_{\rm ou}\!\triangleq
C_{\mathrm{in}}( 1\!+\!\frac{L_{g,1}}{\mu _g} ) ( L_{f,1}\!+\!\frac{L_{g,2}C_{\mathrm{in}}}{2} ). \label{EQ-constant-symbols-cc}
\end{align}

We now introduce   the definition of  $\epsilon$-stationary points \cite{ghadimi2018approximation,lian2017can,yang2022decentralized}, which will be used  in the subsequent analysis to characterize the convergence performance.
\begin{definition}[\textbf{$\epsilon$-stationary point}]
Given a differentiable function  $\Phi:\mathbb{R}^n\rightarrow\mathbb{R}$, we say that $\tilde{x}$ is an $\epsilon$-stationary point  if $\|\nabla\Phi(\tilde{x})\|^2\leqslant\epsilon$.
\end{definition}

In what follows, we make necessary assumptions related to the boundedness of the heterogeneity of the out- and  inner-level functions, which is commonly employed in analyzing the performance of penalty-based methods \cite{nedic2009distributed,lian2017can}.
\begin{definition}[\textbf{Node heterogeneity}]
Given a set of functions $ h_i(x,y): \mathbb{R}^n \times \mathbb{R}^p \rightarrow \mathbb{R}$, $i =1, \cdots, m $, and a constant \( \delta \), we say that the gradients $\{\nabla_x h_i\}$ exhibit
 $\delta ^2$-heterogeneity at any point $x \in \mathbb{R}^{n}$ and  $y \in \mathbb{R}^{p}$
if the condition
$\frac{1}{m}\sum_{i=1}^m{\| \nabla_x h_i\left( x, y\right) - \nabla_x h\left( x, y\right) \|^2} \leqslant \delta ^2
$ holds, where \( h(x, y) = \frac{1}{m} \sum_{i=1}^{m} h_i(x, y) \).
\end{definition}
\begin{assumption}
[\textbf{Bounded  heterogeneity}]  \label{ASS-heterogeneity} There exist some positive constants $b_f^2$ and $b_g^2$ such that: \textbf{(i)} Both $\{\nabla_x f_i\}$ and $\{\nabla_y f_i\}$ have $b_f^2$-heterogeneity at any point   $(x,y^*(x))$ with $x \in \mathbb{R}^n$; \textbf{(ii)} Both $\{\nabla_x g_i\}$ and $\{\nabla_y g_i\}$ have $b_g^2$-heterogeneity at any point $(x,y^*(x))$ with $x \in \mathbb{R}^n$.
\end{assumption}

\begin{remark}[\textbf{Weaker assumption on heterogeneity}]
Assu- mption \ref{ASS-heterogeneity} requires only that the outer-level and inner-level heterogeneity be uniformly bounded at the optimum
$y^*(x)$ for all $x$. This condition is weaker than those in prior works  \cite{chen2023decentralized,yang2022decentralized,kong2024decentralized} which assume that \( \|\nabla_x f_i\| \) is bounded or that the outer-level and inner-level heterogeneity is uniformly bounded for  any point $x \in \mathbb{R}^{n}$ and  $y \in \mathbb{R}^{p}$. Moreover, Assumption \ref{ASS-heterogeneity} focuses exclusively on first-order gradient heterogeneity, unlike earlier studies \cite{niu2023distributed, kong2024decentralized}, which impose additional second-order gradient requirements.
\end{remark}

Assumption 4 is typically satisfied in practical distributed learning scenarios where all nodes share the same objective function structure and their local datasets are drawn from i.i.d. or mildly non-i.i.d. distributions. In such cases, the heterogeneity constants $b_f^2$ and $b_g^2$ quantify the variability across local datasets. This type of bounded heterogeneity assumption is commonly used in the literature on single-level decentralized optimization \cite{lian2017can}.

\subsection{{\ALGNAME} Algorithm}
We are now ready to present the {\ALGNAME} algorithm.
Based on Assumptions \ref{ASS-outer-level} and \ref{ASS-inner-level}, it can be observed that  the subproblem in \eqref{EQ-P-penalty-adptive} is strongly convex in $\{y_i\}$ and strongly concave in $\{z_i\}$ for an appropriately chosen \(\lambda\), with the subproblems for  $\{x_i\}$, $\{y_i\}$, and $\{z_i\}$ exhibiting different levels of complexity. By adopting a loopless update scheme, we expect that the simpler inner-level subproblems can be solved more accurately by using larger step sizes, which, in turn, provides more precise descent directions for the outer-level subproblem. Building on  these insights,
we  propose an
\underline{a}lternating \underline{He}ssian-free \underline{a}pproximate  \underline{d}istributed  algorithm (termed \ALGNAME) in Algorithm \ref{alg:1}
by employing the single-level approximation \eqref{EQ-P-penalty-adptive} along with time-scale separation schemes, which  performs  multiple-timescale and alternating gradient  updates  for the variables  $\{z_i\}$, $\{y_i\}$ and $\{x_i\}$ at each iteration. The compact form of the proposed algorithm  is presented as  follows:
\begin{align}
&\!\!\!\!\!{z^{k+1}=}\mathcal{W}_r z^k-\gamma {\nabla _yG\left( x^k,z^k \right) }, \label{EQ-alg-a}
\\
&\!\!\!\!\!y^{k+1}=\mathcal{W}_r y^k-\beta\left( \nabla _yF\left( x^k,y^k \right) +\lambda \nabla _yG\left( x^k,y^k \right) \right) , \label{EQ-alg-b}  \\
&\!\!\!\!\!x^{k+1}=\mathcal{W}_n x^k \! - \!\alpha
\left( \nabla _xF\left( x^k,y^k \right) \!+\!\lambda \nabla_x Q(x^k,y^k,z^k) \right), \label{EQ-alg-c}
\end{align}
where  $\alpha$, $\beta$, $\gamma$ are the step sizes, with their choices corresponding to different timescales.

\begin{remark}[\textbf{Loopless structure}]
In the proposed {\ALGNAME} algorithm, each node updates \( z_i^{k} \) via a weighted combination of neighboring estimates and a gradient descent step to track \( y^*(\bar{x}^k) \). Then, the value of \( y_i^{k} \) and \( x_i^{k} \) are updated using the gradient of the outer-level function and the penalty term associated with \( \lambda \). The former drives the descent of the outer-level objective, while the latter enforces  the constraint induced by the inner-level subproblem  in  \eqref{EQ-P1}. Note that
  the updates for \( \{z_i^{k}\} \), \( \{y_i^{k}\} \), and \( \{x_i^{k}\} \) do not involve additional computational loops. Unlike existing works \cite{yang2022decentralized, chen2023decentralized, niu2023loopless, dong2023single, niu2023distributed, kong2024decentralized, zhu2024sparkle}, the proposed {\ALGNAME} algorithm is loopless, which, combined with its  Hessian-free property, offers significant computational advantages for large-scale problems.
\end{remark}

\begin{algorithm}[http]
        \caption{\ALGNAME}
	\label{alg:1}
	\begin{algorithmic}[1]
		\STATE \textbf{Require}: Initialize  $\{x_i^0, y_i^0, z_i^0\}$ arbitrarily  and set $\{\alpha, \beta, \gamma, \lambda\}$ according to the selection rules  \eqref{EQ-stepsize-alpha}-\eqref{EQ-stepsize-gamma}.
        \FOR{$k=0, 1, \cdots, K$}
          \STATE \textbf{Update  inner-level variables}:
          $$              z_{i}^{k+1}=\sum \nolimits_{j\in \mathcal{N} _i}{w_{ij}z_{i}^{k}}-\gamma \nabla _yg_i( x_{i}^{k},z_{i}^{k} ),$$
             \STATE \textbf{Update penalized inner-level variables}:
             \begin{align}
              y_{i}^{k+1}=&\!\sum \nolimits_{j\in \mathcal{N} _i}\!{w_{ij}y_{i}^{k}}\nonumber\\
              &- \beta (\nabla _yf_i( x_{i}^{k},y_{i}^{k} ) +\lambda \nabla _yg_i( x_{i}^{k},y_{i}^{k} )),   \nonumber
             \end{align}
            \vspace{-0.5cm}
            \STATE \textbf{Update outer-level variables}:
        \begin{align}
         \!\!\! \!\!x_{i}^{k+1}&=\sum \nolimits_{j\in \mathcal{N} _i}{w_{ij}x_{i}^{k}}  \nonumber \\
          &\!\!-\!\alpha (\nabla _x\!f_i( x_{i}^{k},y_{i}^{k} ) \!+\!\lambda (\!\nabla _xg_i( x_{i}^{k},y_{i}^{k} \!) \!-\!\nabla _xg_i( x_{i}^{k},z_{i}^{k} ) \!)\!).\nonumber
        \end{align}
       \ENDFOR
\end{algorithmic}
\end{algorithm}

\begin{remark}[\textbf{Time-scale separation}] The proposed {\ALGNAME} algorithm performs the updates using  multiple-timescale step sizes \(\alpha\), \(\beta\), and \(\gamma\)  respectively for \(\{x_i^k\}\), \(\{y_i^k\}\), and \(\{z_i^k\}\). In particular, the step sizes \(\beta\) and \(\gamma\) are chosen to decay more slowly than \(\alpha\) in terms of timescale, enabling \(\{y_i^k\}\) and \(\{z_i^k\}\) to converge more quickly and asymptotically to their immediate solutions $y^*(\bar{x}^k;\lambda)$ and $y^*(\bar{x}^k)$ at each iteration given the outer-level variables $\{x_i^k\}$  and
guaranteeing the algorithm's convergence under the loopless structure.
Besides, for the estimates \(\{y_i^k\}\) and \(\{x_i^k\}\), the gradients of the outer-level objective and the penalty term evolve on different timescales, with the penalty term employing a step size that decays more slowly. Such different timescales further ensure the descent of the value of  the outer-level objective while prioritizing the satisfaction of the constraint imposed by the inner-level subproblem at each iteration.
\end{remark}


\subsection{Convergence Results}
In this section, we present the main convergence result
of the proposed {\ALGNAME} algorithm. To this end,  we begin by introducing a potential function as follows:
\begin{align}
\!V^{k}\triangleq & 2\Phi ( \bar{x}^k ) +d_1\| \bar{y}^k\!-\!y^*( \bar{x}^k;\lambda ) \| ^2\!+\!d_2\| \bar{z}^k-y^*( \bar{x}^k ) \| ^2 \nonumber \\
&+\!d_3\frac{1}{m}\| x^{k}\!-\!1_m\otimes \bar{x}^{k} \| ^2 \!+ \!d_4\frac{1}{m}\| y^k-1_m\otimes \bar{y}^k \| ^2\! \nonumber \\
&+\!d_5\frac{1}{m}\| z^{k}\!-\!1_m\otimes \bar{z}^{k} \| ^2, \label{EQ-potential}
\end{align}
where the coefficients $d_1, d_2, d_3, d_4, d_5$ are properly set as: $d_1=
\frac{12U_{\lambda}^2\alpha}{w_{\beta}\beta}
$, $d_2=
\frac{12L_{g,1}^{2}\lambda ^2\alpha}{w_{\gamma}\gamma}$, $
{ d_3=\frac{4p_1\alpha \lambda ^2}{1-\rho}}$, $
d_4={ \frac{4p_3\alpha \lambda ^2}{1-\rho}}
$, and $
d_5={ \frac{4p_3\alpha \lambda ^2}{1-\rho}}
$, with the parameters $p_1$, $p_2$ and $p_3$  given by
\begin{equation}
\begin{aligned} \label{EQ-ppp}
& { p_1=( \frac{96}{u_{\beta}^{2}}+\frac{48L_{g,1}^{2}}{w_{\gamma}^{2}}+24 ) L_{g,1}^{2}}, \\
 &{ p_2=( \frac{96}{u_{\beta}^{2}}+24 ) L_{g,1}^{2}}, \; {p_3=( \frac{48L_{g,1}^2}{w_{\gamma}^{2}}+12 ) L_{g,1}^{2}},
\end{aligned}
\end{equation}
where  $
{u_{\beta}=\frac{\mu_g}{4L_{g,1}}}
$. It should be noted that the introduced potential function primarily consists of four types of dynamics: the objective loss \(\Phi(\bar{x}^k)\), the penalized inner-level errors \(\| \bar{y}^k - y^*(\bar{x}^k; \lambda) \|^2\), the inner-level errors \(\| \bar{z}^k - y^*(\bar{x}^k) \|^2\), and the consensus errors \(\| x^k - 1_m \otimes \bar{x}^k \|^2\), \(\| y^k - 1_m \otimes \bar{y}^k \|^2\), and \(\| z^k - 1_m \otimes \bar{z}^k \|^2\), which provides a comprehensive metric for assessing the convergence behavior  of the sequence $\{x_i^k,y_i^k,z_i^k\}$. Then, we present  the following theorem  to establish   the convergence of the {\ALGNAME} algorithm for solving the nonconvex-strongly-convex  DBO problem \eqref{EQ-P1}.
\begin{theorem}\label{TH-the1}
Consider the sequence $\{x_i^k, y_i^k, z_i^k\}$ generated by Algorithm \ref{alg:1}. Suppose Assumptions  \ref{ASS-outer-level}-\ref{ASS-heterogeneity} hold. If the penalty parameter $\lambda$ satisfies $\lambda>\frac{2L_{f,1}}{\mu_g}$ and the step sizes $\alpha$, $\beta$, $\gamma$ respectively satisfy
\begin{align}
\alpha < \min \Big\{ & \frac{1}{2L},\frac{w_{\gamma}\gamma}{10L_{y^*}L_{g,1}\lambda},\frac{u_{g}^{2}}{160L_{g,1}^{2}}\beta, \label{EQ-stepsize-alpha} \\
&{\frac{1-\rho}{12\lambda \sqrt{p_1}},}\frac{\left( 1-\rho \right) \sqrt{p_3}}{32L_{g,1}\lambda \sqrt{p_1}},\frac{1-\rho}{54L_{g,1}\lambda} \Big\},  \nonumber
\\
\beta <\min \Big\{& \frac{1}{\lambda L_{g,1}},\frac{1}{2(L_{f,1}+\lambda L_{g,1})}+\frac{1}{\lambda \mu_g}, \nonumber\\
&{ \frac{1-\rho}{10\lambda \sqrt{p_2}},}\frac{\left( 1-\rho \right) \sqrt{p_1}}{48L_{g,1}\lambda \sqrt{p_2}},\frac{1-\rho}{32L_{g,1}\lambda} \Big\}, \label{EQ-stepsize-beta}
\\
\gamma <\min \Big\{ &\frac{2}{\mu _g+L_{g,1}},\frac{\mu _g+L_{g,1}}{2\mu _gL_{g,1}},\nonumber\\
&{ \frac{1-\rho}{8\sqrt{p_3}},}\frac{\left( 1-\rho \right) \sqrt{p_1}}{24L_{g,1}\sqrt{p_3}},\frac{1-\rho}{16L_{g,1}} \Big\},  \label{EQ-stepsize-gamma}
\end{align}
then, for any total number of iterations $K$,
we have
\begin{align}
\frac{1}{K}\sum_{k=0}^{K-1}{\left\| \nabla \Phi \left( \bar{x}^k \right) \right\| ^2}
&\leqslant \frac{V^0-V^{K-1}}{\alpha K}  \nonumber\\
&+\underbrace{\mathcal{C}(\lambda, \alpha,\beta, C_{\rm ou}, C_{\rm in},\rho) }_{\rm penalty \; errors} \nonumber\\
&+\underbrace{\mathcal{B}(\lambda,\alpha,\beta,\gamma,b_f, b_g,\rho)}_{\rm heterogeneity \; errors}, \label{EQ-the1}
\end{align}
where $\mathcal{C}(\cdot)$ and $\mathcal{B}(\cdot)$ represent  the penalty and heterogeneity errors respectively, and are given by:
\begin{align}
&\mathcal{C}(\lambda, \alpha,\beta, C_{\rm ou}, C_{\rm in},\rho) \! \label{EQ-p-error} \\
=&\frac{2C_{\rm ou}^2}{\lambda ^2}\!+\!\frac{864C_{\rm in}^{2}L_{g,1}^{2}p_1\lambda ^2\alpha ^2}{\left( 1-\rho \right) ^2}\! +\!\frac{576C_{\rm in}^{2}L_{g,1}^{2}p_2\lambda ^2\beta ^2}{\left( 1-\rho \right) ^2}, \nonumber
 \\
&\!\mathcal{B}(\lambda,\alpha,\beta,\gamma,b_f, b_g,\rho) \label{EQ-b-error}\\
=&\frac{24p_1\lambda ^2\alpha ^2}{\left( 1-\rho \right) ^2}b_{f}^{2} +\frac{24p_3\lambda ^2\gamma ^2}{\left( 1-\rho \right) ^2}b_{g}^{2} \nonumber+\frac{48p_2\lambda ^2\beta ^2}{\left( 1-\rho \right) ^2}\left( b_{f}^{2}+\lambda ^2b_{g}^{2} \right).   \nonumber
\end{align}
\end{theorem}

\begin{proof}
See Section \ref{Sec-poof-main}.
\end{proof}

Theorem \ref{TH-the1} shows$\!$  that  the proposed {\ALGNAME} algorithm can provide a  convergence guarantee\footnote{{The average gradient $\frac{1}{K}\sum_{k=1}^K \|\nabla\Phi(\bar{x}^k)\|^2$ measures the convergence to an $\epsilon$-stationary point. As $K$ increases, if it approaches zero, there exists  an estimate $\bar{x}^t$ with $t \in \{1,\cdots, K\}$ such that $\|\nabla\Phi(\bar{x}^t)\|^2\leqslant\epsilon$.}} for solving the original problem \eqref{EQ-P1}   without requiring  the evaluation of Hessians.
The algorithm's convergence depends on penalty errors \(\mathcal{C}(\lambda, \alpha, \beta, C_{\rm ou}, C_{\rm in}, \rho)\), heterogeneity errors \(\mathcal{B}(\lambda, \alpha, \beta, \gamma, b_f, b_g, \rho)\), initialization states \(V^0\), and step size selection. In particular, the penalty errors for \(C_{\rm ou}\) and \(C_{\rm in}\) in \eqref{EQ-p-error}  are associated with the outer- and inner-level penalty gaps (c.f., Lemma \ref{LE-inner-penalty-error}), respectively, both arising from the penalization of the constraint imposed by the inner-level subproblem. Notably, the penalty errors for \(C_{\rm in}\) is further affected by the penalization of  the consensus constraint, which relies on the   network connectivity and the selection of the step sizes.
The heterogeneity errors in \eqref{EQ-b-error} increase with larger step sizes \(\alpha\), \(\beta\), \(\gamma\) and greater first-order gradient heterogeneity \(b_f\) and \(b_g\). The value of $V^0$ increases with larger values of \(\alpha\) but decreases with larger  \(\beta\) and \(\gamma\). Moreover, these errors exhibit an intricate interaction with the penalty parameter \(\lambda\), as elaborated in the following remark.



\begin{remark}[\textbf{Effect of the penalty parameter}]
We can observe from Theorem \ref{TH-the1} that the penalty parameter \( \lambda \) plays a crucial role in controlling the selection of step sizes $\alpha$, $\beta$ and $\gamma$, and influences both the penalty errors in \eqref{EQ-p-error} and heterogeneity errors in \eqref{EQ-b-error}. Specifically,
it can be seen from the condition \eqref{EQ-stepsize-alpha}-\eqref{EQ-stepsize-gamma} that the selection of both step sizes \( \alpha \) and \( \beta \) is controlled by $\lambda$, with \( \alpha \) being no larger than \( \beta \). Notably, the condition \( \gamma > \frac{10L_{y^*}L_{g,1}}{w_{\gamma}} \alpha \lambda \) from \eqref{EQ-stepsize-alpha} implies that the updates  for the estimates \( \{z_i^k\} \) need to be coordinated with the penalty amplitude of the constraint imposed by  the inner-level subproblem. Regarding the penalty errors, as \( \lambda \) increases, the penalty errors at  the outer level decrease, ultimately approaching zero only as \( \lambda \to \infty \). However, increasing   $\lambda$ and  the step-sizes lead to larger outer- and inner-level heterogeneity errors  (associated with $b_{f}$ and $b_{g}$, respectively), indicating a trade-off between the penalty and heterogeneity errors. Importantly, the inner-level heterogeneity errors  depend on a higher order of $\lambda$, indicating a greater impact of the penalization on the inner-level consensus.
\end{remark}

The result above, combined with the varying levels of difficulty encountered by the variables \(\{x_i\}\), \(\{y_i\}\), and \(\{z_i\}\) in their respective subproblems in \eqref{EQ-P-penalty-adptive}  motivates the detailed selections of step sizes in the following corollary. We now present the convergence rate of the proposed {\ALGNAME} algorithm.

\begin{corollary}\label{COR-rate}
Consider the sequence $\{x_i^k, y_i^k, z_i^k\}$ generated by Algorithm \ref{alg:1}. Suppose Assumptions   \ref{ASS-outer-level}-\ref{ASS-heterogeneity} hold. Given any total number of iterations $K$, if the step sizes $\alpha$, $\beta$ and $\gamma$ are set as $\alpha=\mathcal{O}(\kappa ^{-3}K^{-\frac{2}{3}})$, $\beta=\mathcal{O}(\kappa ^{-1} K^{-\frac{1}{2}})$ and $\gamma=\mathcal{O}(K^{-\frac{1}{3}})$, respectively, and $\lambda=\mathcal{O}(\kappa K^{\frac{1}{6}})$ such that the condition \eqref{EQ-stepsize-alpha}-\eqref{EQ-stepsize-gamma} holds, then we have
\begin{align}
&\frac{1}{K}\sum _{k=0}^{K-1}{\| \nabla \Phi ( \bar{x}^k ) \| ^2} \label{EQ-cor1-1} \\
 &~~~\leqslant
\mathcal{O} \big( \frac{\kappa ^4}{K^{1/3}}+\frac{1}{(1-\rho)^2}(\frac{{\kappa ^2b_{f}^{2}}}{K^{2/3}}+\frac{\kappa ^4b_{g}^{2}}{K^{1/3}}) \big) , \nonumber
\end{align}
where $\kappa=\frac{\max\{C_{f,y}, L_{f,1},L_{g,1},L_{g,2}\}}{\mu _g}$ is the condition number.
\end{corollary}
\begin{proof}
See Section \ref{Sec-poof-main}.
\end{proof}


Corollary \ref{COR-rate} shows
that the {\ALGNAME} algorithm achieves a  convergence rate of \(\mathcal{O}(\frac{\kappa^4}{(1-\rho)^2 K^{1/3}} )\)  under time-scale separation for the update of $\{x_i^k\}$, $\{y_i^k\}$ and $\{z_i^k\}$. The result  clearly reveals the detailed influence of the inner-level heterogeneity \(b_f^2\), outer-level heterogeneity \(b_g^2\), network connectivity \(\rho\), and  condition number \(\kappa\) on the convergence. The dependence of the convergence rate on $\rho$ in Corollary~\ref{COR-rate} is consistent with that observed in    the DGD method \cite{nedic2017achieving, li2018distributed, lian2017can} for distributed single-level optimization. By choosing step sizes that depend on $\rho$, such as \(\alpha = \mathcal{O}(1 - \rho)\), the dependence of the  rate on network connectivity can be reduced to \(\mathcal{O}(\frac{1}{1 - \rho})\).
It should be noted that, since only first-order gradients are used, the convergence rate of the {\ALGNAME} algorithm is slower  compared to those based on Hessian evaluations for DBO which achieve a rate of \(\mathcal{O}(\frac{1}{K^{2/3}})\) with the DGD method \cite{kong2024decentralized,niu2023distributed} and a faster rate of \(\mathcal{O}(\frac{1}{K})\) with the GT method \cite{dong2023single,niu2023distributed}. Due to the impact of the heterogeneity and the elimination of the additional computation loop, the obtained rate is slower than those in centralized Hessian-free approaches \cite{liu2022bome, kwon2023fully, yao2024constrained}, which achieve a rate of \(\mathcal{O}(\frac{1}{K^{2/3}})\) under additional computation loops. In particular, as shown in Theorem \ref{TH-the1}, the additional heterogeneity errors are introduced in DBO, with the inner-level heterogeneity errors being strongly affected by the penalty parameter \( \lambda \). To ensure convergence of the inner-level heterogeneity errors, a smaller \( \lambda \) is required, which slows down the convergence  of the penalty errors, resulting in slower overall convergence.



\begin{remark}[\textbf{Computational complexity}] Note that the proposed {\ALGNAME} algorithm is able to achieve a computational complexity\footnote{The computational complexities are measured by the number of outer-level and inner-level gradients and their dimensions.} of \(\mathcal{O}(n + r)\) per iteration \(k\) and a total computational complexity of \(\mathcal{O}((n + r) \epsilon^{-3})\) to reach an \(\epsilon\)-stationary point, due to its Hessian-free and loopless structure.  In comparison, existing methods \cite{yang2022decentralized, chen2023decentralized} incur a computational complexity of \(\mathcal{O}((n^2 + nr)\log{k})\) per iteration \(k\), leading to a total complexity of \(\mathcal{O}((n^2 + nr) \epsilon^{-3/2} \log(\epsilon^{-1}))\). Similarly, the methods in \cite{dong2023single, niu2023distributed} require \(\mathcal{O}((n^2 + nr) \epsilon^{-3/2})\), while those in \cite{dong2023single, niu2023distributed, zhu2024sparkle} need \(\mathcal{O}((n^2 + nr) \epsilon^{-1})\). We see that the {\ALGNAME} algorithm enjoys  a significant advantage in computational complexity per iteration. Notably, for large-scale problems, high solution accuracy is not essential (i.e., $\epsilon$ being relatively large), such as in machine learning, the proposed algorithm offers  more  efficient overall computational complexity compared to existing methods \cite{yang2022decentralized, chen2023decentralized, niu2023loopless, dong2023single, niu2023distributed, kong2024decentralized, zhu2024sparkle}, although it exhibits a slower convergence rate.
\end{remark}

\begin{remark}[\textbf{Heterogeneity analysis}] \label{RE-revised-heter} It is shown in
Corollary \ref{COR-rate} that the error term related to the heterogeneity of the inner-level functions, \(b_g^2\), converges more slowly than the error term related to \(b_f^2\), and thus has a dominant impact on the convergence performance of the heterogeneity errors when $K$ is relatively large. This result underscores the significant role of   the inner-level heterogeneity when applying the penalty method to  DBO. It is noted that the heterogeneity analysis in existing studies \cite{ niu2023distributed, kong2024decentralized}  relies on the condition on the second-order  heterogeneity. In contrast,  our heterogeneity analysis is based solely on  first-order heterogeneity  \(b_f^2\) and \(b_g^2\) (c.f., Assumption \ref{ASS-heterogeneity}) similar to that used in single-level optimization \cite{lian2017can}, which allows us to provide a tighter analysis for the impact of the heterogeneity.
To our best knowledge, our work
 is the first to  establish such a heterogeneity result for DBO, which bridges the gap between DBO and distributed  single-level optimization in terms of heterogeneity.
\end{remark}

\subsection{Applications to Distributed Minimax Problems}
In this section, we extend our results to  distributed minimax problems, which have important practical applications in scenarios such as distributed adversarial attacks, fair beamforming, and the training of generative adversarial networks \cite{razaviyayn2020nonconvex}.
In particular, when the local inner-level objective exactly opposes  the local outer-level objective, i.e., \(g_i(x, y) = -f_i(x, y)\) holds for all \(i \in \mathcal{V}\), the problem \eqref{EQ-P1} will reduce to a special case: \(\min_x \Phi(x) = \min_x (1/m)\sum_{i=1}^{m} f_i(x, y^*(x))\) with \(y^*(x) = \arg\max_y (1/m)\sum_{i=1}^{m} f_i(x, y)\). This case is equivalent to the following distributed minimax problem when $f_i(x,y)$ is strongly concave in $y$:
\begin{equation}\label{EQ-minmax}
\begin{aligned}
&\underset{x\in \mathbb{R} ^n}{\min} \underset{y \in \mathbb{R} ^r}{\max}\frac{1}{m}{\sum_{i=1}^m{f_i\left( x,y\right)}}.
\end{aligned}
\end{equation}
In this case, due to the strong concavity of \(f_i(x, y)\) in \(y\), we have \(\nabla_y f(x, y^*(x)) = 0\), which simplifies the expression of the hypergradient as \(\nabla \Phi(x) = \nabla_x f(x, y^*(x))\).
Under Assumptions \ref{ASS-outer-level} and \ref{ASS-inner-level}, we can refine the result of Proposition \ref{PRO-PHI} for \(\nabla \Phi(x)\), showing that \(\nabla \Phi(x)\) is \(L\)-Lipschitz continuous with \(L = L_{f,1}+L_{f,1}L_{y^*}\). Moreover, in the distributed minimax case, the updates of the variables $\{x_i^k\}$ and $\{y_i^k\}$ in Algorithm \ref{alg:1}  can be simplified as:
   \begin{align}
         \!\!\! \!\!z_{i}^{k+1}=&\sum \nolimits_{j\in \mathcal{N} _i}{w_{ij}z_{i}^{k}} +\!\gamma \nabla _y f_i( x_{i}^{k},z_{i}^{k} \!), \label{EQ-minmax-z}\\
         \!\!\! \!\!y_{i}^{k+1}=&\sum \nolimits_{j\in \mathcal{N} _i}{w_{ij}y_{i}^{k}} -\!\beta (1-\lambda)\nabla _y f_i( x_{i}^{k},y_{i}^{k} \!), \label{EQ-minmax-y}\\
         \!\!\! \!\!x_{i}^{k+1}=&\sum \nolimits_{j\in \mathcal{N} _i}{w_{ij}x_{i}^{k}}  \label{EQ-minmax-x}\\
          &-\!\alpha ((1-\lambda)\nabla _x f_i( x_{i}^{k},y_{i}^{k} ) +\lambda \nabla _x f_i ( x_{i}^{k},z_{i}^{k} )),\nonumber
        \end{align}
whereby we can see that   the variables $\{y_i^k\}$ and $\{z_i^k\}$ in the {\ALGNAME} algorithm now aim to estimate the same inner-level solution.  Notably, when \(\lambda = 1\) and the step sizes \(\beta\) and \(\gamma\) are set equal, with the same initial values for \(y_i^0\) and \(z_i^0\), the recursions in \eqref{EQ-minmax-x}-\eqref{EQ-minmax-y} reduce to the standard distributed gradient descent-ascent updates for the distributed minimax problem, which further demonstrates the versatility  of the proposed {\ALGNAME} algorithm.
The following corollary provides  convergence results of the {\ALGNAME} algorithm  for the distributed minimax case.

\begin{corollary}[\textbf{Distributed minimax problems}]\label{COR-cor2}
Consider the case where \(g_i(x, y) = -f_i(x, y)\) for all \(i \in \mathcal{V}\), and the sequence \(\{x_i^k, y_i^k, z_i^k\}\) generated by Algorithm \ref{alg:1}.  Suppose Assumptions  \ref{ASS-outer-level}, \ref{ASS-network}, \ref{ASS-heterogeneity} and  Assumptions \ref{ASS-inner-level}(i), \ref{ASS-inner-level}(ii) hold. Let \(L = L_{f,1}+L_{f,1}L_{y^*}\). If  $\lambda
\geqslant 1$, then we have that  the problem \eqref{EQ-minmax} becomes equivalent to its penalty formulation $\min_{x,y}\max_z p(x,y,z;\lambda)$ and the inner-level penalty gap  reduces to zero, i.e.,  ${\left\| y^*\left( x \right) -y^*\left( x;\lambda \right) \right\| ^2=0}$. Further, if the step sizes $\alpha$,  $\beta$, and $\gamma$  are set as $\alpha=\mathcal{O}(\kappa^{-2}K^{-{\frac{1}{3}}})$, $\beta=\mathcal{O}(K^{-{\frac{1}{3}}})$, and $\gamma=\mathcal{O}(K^{-{\frac{1}{3}}})$, with  \(\lambda\alpha = \mathcal{O}(\kappa^{-2}K^{-{\frac{1}{3}}})\) and  \(\lambda\beta = \mathcal{O}(K^{-{\frac{1}{3}}})\)\footnote{{Here, the conditions $\alpha=\mathcal{O}(\kappa^{-2}K^{-{\frac{1}{3}}})$ and \(\lambda\alpha = \mathcal{O}(\kappa^{-2}K^{-{\frac{1}{3}}})\) require that the parameter $\lambda$ is selected independently of both the number of iterations and the condition number.}},  such that the condition \eqref{EQ-stepsize-alpha}-\eqref{EQ-stepsize-gamma} are satisfied, then  we have
\begin{equation}
\begin{aligned}\label{EQ-COR2-2}
 \!\!\!\! \frac{1}{K}\sum_{k=0}^{K-1}{\| \nabla \Phi \left( \bar{x}^k \right) \| ^2}
\! \leqslant
 \!
\mathcal{O}\big(\! \frac{\kappa^2}{ K^{\frac{2}{3}}  }\!\!+\! \!\frac{\kappa^2 b_{f}^{2}}{\left( 1-\rho \right) ^2  K^{\frac{2}{3}} } \big).
\end{aligned}
\end{equation}
\end{corollary}
\begin{proof}
See Section \ref{Sec-poof-main}.
\end{proof}

\begin{remark}[\textbf{Tight convergence analysis}]
From Corollary \ref{COR-cor2}, we observe that   if  the penalty parameter $\lambda$ satisfies \(\lambda \geqslant 1\), the approximation in \eqref{EQ-P-penalty-adptive} can be treated as an exact penalty for the constraint \(y^*(x) = \arg\max_y \frac{1}{m} \sum_{i=1}^{m} f_i(x, y)\) in the distributed minimax case.
Besides, Corollary  \ref{COR-cor2} demonstrates that the proposed {\ALGNAME} algorithm  can achieve a convergence rate of $\mathcal{O}({\frac{\kappa^2}{(1-\rho)^2K^{\frac{2}{3}}}})$ for distributed minimax problems.  When the step sizes $\alpha$,  $\beta$, and $\gamma$ are chosen as the constants independent of $K$, i.e., $\alpha=\mathcal{O}(\kappa^{-2})$, $\beta=\mathcal{O}(1)$, and $\gamma=\mathcal{O}(1)$ under the condition \eqref{EQ-stepsize-alpha}-\eqref{EQ-stepsize-gamma}, we have
$ \! \frac{1}{K}\sum_{k=0}^{K-1}{\| \nabla \Phi \left( \bar{x}^k \right) \| ^2}
 \leqslant
\mathcal{O}\big(\! \frac{\kappa^2}{ K  }+\frac{\kappa^2 b_{f}^{2}}{\left( 1-\rho \right) ^2} \big)
$, which implies a faster rate of $\mathcal{O}(\frac{\kappa ^2}{K})$ to  a neighborhood of the stationary point. This result is consistent with that of the state-of-the-art work \cite{lin2020gradient},  further highlighting the tightness of our convergence analysis.
\end{remark}

\section{Convergence Analysis}
In this section, we present the convergence analysis of the main results for the {\ALGNAME} algorithm. Our proof aims to establish the evolution of the potential function in \eqref{EQ-potential} by analyzing the contraction properties of the average system, which are related to the penalized inner-level errors, the inner-level errors, and the dynamics of the consensus errors. This analysis enables us to further develop the results in Theorem \ref{TH-the1} and Corollaries \ref{COR-rate} and \ref{COR-cor2}. To facilitate the subsequent convergence analysis, we first provide some technical preliminaries.

\subsection{Technical Preliminaries}
For brevity, we introduce some compact notations, including the average  operator $J_n \triangleq \frac{1_m^{\rm{T}} \otimes I_n}{m}$, the consensus operator $\mathcal{J}_n \triangleq \frac{1_m 1_m^{\rm{T}}}{m} \otimes I_n$, along with the  gradient expressions  of the {\ALGNAME} algorithm:
\begin{align}
&H_{z}^k\triangleq {\nabla _yG\left( x^k,z^k \right) },   \nonumber\\
&H_{y}^k\triangleq  \nabla _yF\left( x^k,y^k \right) +\lambda \nabla _yG\left( x^k,y^k \right),  \label{EQ-HHH} \\
&H_{x}^k\triangleq  \nabla _xF\left( x^k,y^k \right) \!+\!\lambda \nabla_x Q(x^k,y^k,z^k), \nonumber
\end{align}
and  we use bar notation to
denote averages in later analysis, e.g., $\bar{H}_x^k=J_n H_x^k$. We then present the following two propositions concerning the Lipschitz continuity of the hypergradient, the gradient of the penalty function, and their corresponding inner-level solutions.
\begin{proposition}[\cite{ghadimi2018approximation}]\label{PRO-PHI}
Suppose Assumptions \ref{ASS-outer-level} and \ref{ASS-inner-level} hold. Then,  $\nabla \Phi(x)$ is   $L$-Lipschitz-continuous and  $y^*(x)$ is $L_{y^*}$-Lipschitz-continuous.
\end{proposition}

\begin{proposition}\label{PRO-penalty}
Suppose Assumptions \ref{ASS-outer-level} and \ref{ASS-inner-level} hold.
If $\lambda>\frac{2L_{f,1}}{\mu _g}$, then
$\nabla_y p(x,y,z;\lambda)$ is $L_{\lambda}$-Lipschitz-continuous and $y^*(x;\lambda)$ is $L_{y^*,\lambda}$-Lipschitz-continuous.
\end{proposition}

The proof of Proposition \ref{PRO-penalty} follows directly from the Lipschitz continuity of \( \nabla f_i \) and \( \nabla g_i \), along with the strong convexity of the penalty function \( p(x, y, z; \lambda) \) in \( y \). With these propositions and notations established, we now  prove the main results for the {\ALGNAME} algorithm.

\subsection{Supporting  Lemmas}
We aim to prove the convergence of the {\ALGNAME} algorithm in terms of the first-order stationary point metric for the average system, i.e., \(\nabla \Phi(\bar{x}^k)\). To this end, we first characterize how the objective function \(\Phi\) decreases along the sequence of averaged iterates \(\{\bar{x}^k\}\).
\begin{lemma}[\textbf{Descent lemma}]\label{LE-descent}
Consider the sequence $\{x_i^k, y_i^k, z_i^k\}$ generated by Algorithm \ref{alg:1}. Suppose Assumptions \ref{ASS-outer-level}-\ref{ASS-network} hold. If the step-size satisfies
$
\alpha\leqslant\frac{1}{2L}
$, then we have
\begin{align}
  &\Phi ( \bar{x}^{k+1} ) \\ \leqslant& \Phi ( \bar{x}^k ) \! +\! \frac{\alpha}{2}{\| \nabla \Phi ( \bar{x}^k ) \!- \! \bar{H}_{x}^{k} \| ^2} \!- \!\frac{\alpha}{2}\| \nabla \Phi ( \bar{x}^k ) \| ^2 \!-\!\frac{\alpha}{4}\| \bar{H}_{x}^{k} \| ^2. \nonumber
\end{align}
\end{lemma}
\begin{proof}
See  Appendix \ref{S-p1}.
\end{proof}


The above lemma demonstrates that the decrease of the objective function \(\Phi\) is dominated by the gradient approximation errors \(\|\nabla \Phi (\bar{x}^k) - \bar{H}_x^k\|^2\). To analyze these errors, we first examine the penalty gap introduced by penalizing the constraint associated with the inner-level subproblem in the following lemma.

\begin{lemma}[\textbf{Penalty gap at inner and outer levels}] \label{LE-inner-penalty-error}
Supp-   \\ose Assumptions \ref{ASS-outer-level} and \ref{ASS-inner-level} hold.  Recall that {{$\mu _{\lambda}=\frac{\lambda \mu_g}{2}$}} and $y^*\left( x;\lambda \right) \in \arg\min _y\max _zp\left( x,y,z;\lambda \right)$.   Let $p^*\left( x;\lambda \right)= p(x,y^*(x;\lambda),y^*(x);\lambda)$.
If $\lambda >\frac{2L_{f,1}}{\mu _g}$,  the auxiliary penalty function $p(x,y,z;\lambda)$  is  $
 \mu _{\lambda}$-strongly convex in $y$  for any $x \in \mathbb{R}^n$ and $z \in \mathbb{R}^r$, and the solution $y^*\left( x;\lambda \right)$ to the problem $\min _y\max _zp\left( x,y,z;\lambda \right)$ is unique. Then, the penalty gap between  $y^*(x)$ and  $y^*(x;\lambda)$ is bounded by
\begin{equation}
{\left\| y^*\left( x \right) -y^*\left( x;\lambda \right) \right\| ^2\leqslant \frac{C_{\mathrm{in}}^{2}}{\lambda ^2}},
\end{equation}
and the penalty gap between $\nabla \Phi(x)$   and   $\nabla p^*\left( x;\lambda \right) $ is bounded by
\begin{equation} \label{EQ-approx-error}
\begin{aligned}
\left\| \nabla \Phi \left( x \right) -\nabla p^*\left( x;\lambda \right) \right\|^2 \leqslant{ \frac{C_{\rm ou}^2}{\lambda^2}}.
\end{aligned}
\end{equation}
\end{lemma}

\begin{proof}
See  Appendix \ref{S-p1-0}.
\end{proof}


It is noted that Lemma \ref{LE-inner-penalty-error} is derived under the weaker condition in Assumption \ref{ASS-outer-level}, which only requires the boundedness of \(\|\nabla_y g_i\|\) at the optimum \(y^*(x)\), instead of requiring this condition to be hold for all \(y\).


\begin{lemma}[\textbf{Gradient approximation  errors}]\label{LE-hypergradient-approx}
Consider the sequence $\{x_i^k, y_i^k, z_i^k\}$ generated by Algorithm \ref{alg:1}. Suppose Assumptions \ref{ASS-outer-level}-\ref{ASS-network} hold. Recall that $U_{\lambda}^2 = L_{f,1}^{2}+\lambda ^2L_{g,1}^{2}$. Then,  we have
\begin{align}
&{\| \nabla \Phi ( \bar{x}^k ) -\bar{H}_{x}^{k} \| ^2}   \label{EQ-phi-hx-bar}
\\
\leqslant &
\frac{2C_{\rm ou}^2}{\lambda ^2}+12U_{\lambda}^2 \| \bar{y}^k-y^*( \bar{x}^k;\lambda  ) \| ^2 \nonumber \\
&+12L_{g,1}^{2}\lambda ^2\| \bar{z}^k-y^*( \bar{x}^k ) \| ^2
+12U_{\lambda}^2\frac{1}{m}\| x^k-1_m\otimes \bar{x}^k \| ^2 \nonumber
\\
&+12U_{\lambda}^2 \frac{1}{m}\| y^k\!-\!1_m\otimes \bar{y}^k \| ^2
\!+12L_{g,1}^{2}\lambda ^2\frac{1}{m}\| z^k\!-\!1_m\otimes \bar{z}^k \| ^2. \nonumber
\end{align}
\end{lemma}
\begin{proof}
See Appendix \ref{S-p2}.
\end{proof}


\begin{lemma}[\textbf{Inner-level  errors}]\label{LE-inner-errors}
Consider the sequence $\{x_i^k, y_i^k, z_i^k\}$ generated by Algorithm \ref{alg:1}. Suppose Assumptions \ref{ASS-outer-level}-\ref{ASS-network} hold. Recall that $
{w_{\gamma}=\frac{\mu _gL_{g,1}}{2(\mu _g+L_{g,1})}}
$. If the step-size satisfies
\begin{equation}\label{EQ-condition-gamma}
\gamma \leqslant \min \{ \frac{2}{\mu _g+L_{g,1}},\frac{\mu _g+L_{g,1}}{2\mu _gL_{g,1}} \},
\end{equation}
then we have

\begin{align}
&\| \bar{z}^{k+1}-y^*( \bar{x}^{k+1} ) \| ^2   \label{EQ-inner-errors} \\
\leqslant&
( 1-\frac{\mu _gL_{g,1}}{\mu _g+L_{g,1}}\gamma ) \| \bar{z}^k-y^*( \bar{x}^k ) \| ^2+\frac{2L_{y^*}^{2}}{\gamma w_{\gamma}}\alpha ^2\| \bar{H}_{x}^{k} \| ^2 \nonumber
\\
&+\frac{4L_{g,1}^{2}}{w_{\gamma}}\gamma\frac{1}{m} ( \| x^k-1_m\otimes \bar{x}^k \| ^2+\| z^k-1_m\otimes \bar{z}^k \| ^2 ) \nonumber
.
\end{align}
\end{lemma}
\begin{proof}
See Appendix \ref{S-p3}.
\end{proof}


\begin{lemma}[\textbf{Penalized inner-level  errors}]\label{LE-inner-penalty}
Consider the sequence $\{x_i^k, y_i^k, z_i^k\}$ generated by Algorithm \ref{alg:1}.  Suppose Assumptions \ref{ASS-outer-level}-\ref{ASS-network} hold.
Recall that $
{ w_{\beta}=\frac{\mu _{\lambda}L_{\lambda}}{2(\mu _{\lambda}+L_{\lambda})}}$.
If the step-size satisfies
\begin{equation}\label{EQ-condition-beta}
\begin{aligned}
\beta \leqslant&
\min \{ \frac{1}{\lambda L_{g,1}},{ \frac{1}{2(L_{f,1}+\lambda L_{g,1})}+\frac{1}{\lambda \mu_g}  } \},
\end{aligned}
\end{equation}
then we have
\begin{align}
&\| \bar{y}^{k+1}-y^*( \bar{x}^{k+1};\lambda  ) \| ^2 \label{EQ-inner-penalty}
\\
\leqslant &( 1-\frac{\mu _{\lambda}L_{\lambda}}{\mu _{\lambda}+L_{\lambda}}\beta ) \| \bar{y}^k-y^*( \bar{x}^k;\lambda ) \| ^2 +\frac{2L_{y*,\lambda}^{2}}{\beta w_{\beta}}\alpha ^2\| \bar{H}_{x}^{k} \| ^2 \nonumber
\\
& \!\!+\frac{8U_{\lambda}^2}{w_{\beta}}\!\beta \frac{1}{m}(  \left\| x^k\!-\!1_m \! \otimes \! \bar{x}^k\! \right\| ^2\!+\!\|y^k\!-\!1_m\!\otimes\! \bar{y}^k \!\| ^2\! ) . \nonumber
\end{align}
\end{lemma}
\begin{proof}
See Appendix \ref{S-p4}.
\end{proof}

As shown in Lemmas \ref{LE-inner-errors} and \ref{LE-inner-penalty}, the inner-level errors and penalized inner-level errors, respectively estimated at the averaged iterates \(\bar{z}^k\) and \(\bar{y}^k\), will gradually converge to zero as the sequence of the averaged iterates \(\bar{x}^k\) progresses and the consensus errors are efficiently controlled. In particular, we can also see that the penalty parameter \(\lambda\) will affect the choice of the step-size \(\beta\) for the variable \(\{y_i^k\}\).

\begin{lemma}[\textbf{Consensus errors}]\label{LE-consensus errors}
Consider the sequence $\{x_i^k, y_i^k, z_i^k\}$ generated by Algorithm \ref{alg:1}. Recall that \(\rho =\| W - \frac{1_m 1_m^{\rm{T}}}{m} \|^2 \in [0, 1)\). Suppose Assumptions  \ref{ASS-outer-level}-\ref{ASS-heterogeneity} hold.
Then, we have
\begin{align}
&\| x^{k+1}-1_m\otimes \bar{x}^{k+1} \| ^2
\leqslant ( 1\!-\!\frac{1\!-\!\rho}{2} ) \!\| x^k-1_m\otimes \bar{x}^k \| ^2  \label{EQ-consensus-x}
 \\
&\;\;\;\;\;\;\;\;\;\;\;\;\;\;\;\;\;\;\;\;+\frac{108m\alpha ^2C_{\rm in}^{2}U_{\lambda}^{2}}{( 1-\rho ) \lambda ^2}\! +\frac{6m\alpha^2b_{f}^{2}}{1-\rho}+E_x^c+E_x^o,
\nonumber
\end{align}
\begin{align}
&\!\!\| y^{k+1}\!-1_m\otimes \bar{y}^{k+1} \| ^2 \!
\leqslant
( 1-\frac{1-\rho}{2} ) \left\| y^k\!-\!1_m\otimes \bar{y}^k \right\| ^2\label{EQ-consensus-y}
 \\
&\!\;\;\;\;\;\;\;\;\;\; +\frac{72m\beta ^2C_{\rm in}^2{U_{\lambda}^2}}{(1-\rho)\lambda^2}+\frac{12m\beta ^2}{1-\rho}\left( b_{f}^{2}+\lambda ^2b_{g}^{2} \right)+E_y^c+E_y^o, \nonumber
\end{align}
and
\begin{align}
&\!\!\| z^{k+1}-1_m\otimes \bar{z}^{k+1} \| ^2  \label{EQ-consensus-z}     \\
\leqslant &
( 1-\frac{1-\rho}{2} ) \| z^k-1_m\otimes \bar{z}^k \| ^2+E_z^c+E_z^o+\frac{6m\gamma ^2b_{g}^{2}}{1-\rho}. \nonumber
\end{align}
where the terms
$E_x^c=\frac{108\alpha ^2U_{\lambda}^{2}}{1-\rho}\| y^k\! - \! 1_m\otimes \bar{y}^k \| ^2
+\frac{108\alpha ^2U_{\lambda}^{2}}{1-\rho}\| x^k \!- \!1_m\!\otimes\! \bar{x}^k \| ^2\!+\!\frac{108\alpha ^2\lambda ^2L_{g,1}^{2}}{1-\rho}\| z^k\!-\!1_m\otimes  \bar{z}^k  \| ^2$, $E_y^c\triangleq\frac{72\beta ^2{U_{\lambda}^2}}{1-\rho}\| x^k \!- \!1_m\otimes \bar{x}^k \| ^2+\frac{72\beta ^2{U_{\lambda}^2}}{1-\rho}\| y^k\!-\!1_m\otimes \bar{y}^k \| ^2$, and $E_z^c\triangleq\frac{24L_{g,1}^{2}\gamma ^2}{1-\rho}\| x^k-1_m\otimes \bar{x}^k \| ^2
+\frac{24L_{g,1}^{2}\gamma ^2}{1-\rho}\| z^k-1_m\otimes \bar{z}^k \| ^2$ denote the auxiliary variables related to consensus errors;  the terms $E_x^o=\!\frac{108m\alpha ^2U_{\lambda}^{2}}{1-\rho}\| \bar{y}^k \!- \!y^*( \bar{x}^k;\!\lambda ) \| ^2+\!\frac{108m\alpha ^2\lambda ^2L_{g,1}^{2}}{1-\rho}\| \bar{z}^k\!-\! y^*(\bar{x}^k)  \| ^2$, $E_y^o\triangleq\frac{72m\beta ^2{U_{\lambda}^2}}{1-\rho}\| \bar{y}^k\!-\!y^*( \bar{x}^k;\lambda ) \| ^2$, and $E_z^o\triangleq\frac{24mL_{g,1}^{2}\gamma ^2}{1-\rho}\| \bar{z}^k-y^*( \bar{x}^k ) \| ^2$ denote the auxiliary variables related to the inner-level optimization errors.

\end{lemma}
\begin{proof}
See Appendix \ref{S-p5}.
\end{proof}

Lemma \ref{LE-consensus errors} explicitly shows how the evolution of the consensus errors for $\{x_i^k\}$, $\{y_i^k\}$, and $\{z_i^k\}$ is influenced by the penalty parameter $\lambda$, node heterogeneity $b_f^2$, $b_g^2$, network connectivity $\rho$, and inner-level optimization errors.
In particular, although the gradient used for updating  the outer-level variables in \eqref{EQ-alg-c} is associated with the penalty parameter \(\lambda\), we mitigate  the influence of the penalty parameter $\lambda$  on the errors induced by  the inner-level heterogeneity $b_g^2$ by exploring the potential relationship between the estimates \(\{y_i^k\}\) and \(\{z_i^k\}\). This results in a tighter bound term of order $\mathcal{O}(\alpha^2 b_f^2)$ for the errors induced by heterogeneity in \eqref{EQ-consensus-x}, instead of the order $\mathcal{O}(\alpha^2 (b_f^2 + \lambda^2 b_g^2))$.
Additionally, by leveraging the relationship between \( y^*(\bar{x}^k) \) and \( y^*(\bar{x}^k; \lambda) \), and introducing the inner-level errors and penalized inner-level error terms, we  tightly characterize the evaluation of the consensus errors associated with \( \{x_i^k\} \), \( \{y_i^k\} \), and \( \{z_i^k\} \) under the weaker Assumption \ref{ASS-heterogeneity} on node heterogeneity, which relies
solely on first-order heterogeneity at the optimum for both the outer- and inner-level functions. This technical analysis marks a significant departure from existing works on DBO \cite{niu2023distributed, kong2024decentralized, dong2023single, zhu2024sparkle},  which either lack an explicit heterogeneity analysis or rely on the additional condition of second-order heterogeneity.

\subsection{Proofs of The Main Results} \label{Sec-poof-main}
\colorbf{\textbf{1) Proof of Theorem  \ref{TH-the1}}}. The proof is based on a key descent lemma and the evolution of the gradient approximation errors for the non-convex outer-level objective \(\Phi(\bar{x}^k)\) (c.f., Lemmas \ref{LE-descent}, \ref{LE-inner-penalty-error}, \ref{LE-hypergradient-approx}), the contraction properties of the  errors  related to the inner level (c.f., Lemmas \ref{LE-inner-errors}, \ref{LE-inner-penalty}), and the dynamics of the consensus errors (c.f., Lemma \ref{LE-consensus errors}). With these results, we establish a tight descent property for the carefully designed potential function defined in \eqref{EQ-potential}  by properly selecting the coefficients \(d_1\) to \(d_5\) and applying small-gain-like techniques, ensuring that the overall convergence behavior is well  characterized.

First of all, invoking Lemmas \ref{LE-descent}-\ref{LE-inner-penalty},
and combining them with the definitions of the coefficients \(d_1\) and \(d_2\) in \eqref{EQ-potential}, we can establish the descent property of the outer-level objective and the  errors related to the inner level as follows:
\begin{align}
&\!\!2\Phi ( \bar{x}^{k+1} ) \! + \!d_1\| \bar{y}^{k+1}\!-\!y^*( \bar{x}^{k+1};\lambda ) \| ^2+ \!d_2\| \bar{z}^{k+1}\!-\!y^*( \bar{x}^{k+1} ) \| ^2
\nonumber\\
\leqslant &2\Phi ( \bar{x}^k )+d_1\| \bar{y}^k-y^*( \bar{x}^k;\lambda ) \| ^2 +d_2\| \bar{z}^k-y^*( \bar{x}^k ) \| ^2
\nonumber\\
&-\alpha \| \nabla \Phi ( \bar{x}^k ) \| ^2+\frac{2C^2_{\rm ou}}{\lambda ^2}\alpha
\nonumber\\
&-( \frac{1}{2}-\frac{24L_{y^*}^{2}L_{g,1}^{2}\lambda^2}{w_{\gamma}^{2}\gamma ^2}\alpha ^2-\frac{24U_{\lambda}^{2}L_{y*,\lambda}^{2}}{w_{\beta}^{2}\beta ^2}\alpha ^2 ) \alpha \| \bar{H}_{x}^{k} \| ^2
\nonumber\\
&   -12U_{\lambda}^{2}\alpha\| \bar{y}^k-y^*( \bar{x}^k;\lambda ) \| ^2-12L_{g,1}^{2}\lambda ^2\alpha\| \bar{z}^k-y^*\left( \bar{x}^k \right) \| ^2
\nonumber\\
&+\frac{D_1\alpha}{m}\| x^k-1_m\otimes \bar{x}^k \| ^2+\frac{D_2\alpha}{m}\| y^k-1_m\otimes \bar{y}^k \| ^2 \nonumber \\
&+\frac{D_3\alpha}{m}\| z^k-1_m\otimes \bar{z}^k \| ^2   \label{EQ-d0-d2-V},
\end{align}
where $D_1=\frac{96U_{\lambda}^{4}} {w_{\beta}^{2}}+\frac{48L_{g,1}^{4}\lambda ^2}{w_{\gamma}^{2}}+12U_{\lambda}^{2},
D_2=\frac{96U_{\lambda}^{4}}{w_{\beta}^{2}}+12U_{\lambda}^{2},
D_3=\frac{48L_{g,1}^{4}\lambda ^2}{w_{\gamma}^{2}}+12L_{g,1}^{2}\lambda ^2$.
In what follows, we aim to establish bounds for the terms $D_1$, $D_2$ and $D_3$. By the condition $
{ \lambda \geqslant }\frac{{2L_{f,1}}}{{\mu _g}}
$, it follows that $
L_{f,1}^{2}\leqslant \frac{\lambda ^2\mu _{g}^{2}}{4}\leqslant \lambda ^2L_{g,1}^{2}
$ and $U_{\lambda}^2\leqslant
2L_{g,1}^{2}\lambda ^2$. Furthermore, by the definition of $w_{\beta}$ in \eqref{EQ-constant-symbols}, we have $
\frac{L_{g,1}\lambda}{w_{\beta}}=\frac{2L_{g,1}}{\mu_g}( \frac{\mu _{\lambda}}{L_{\lambda}}+1 ) \leqslant \frac{1}{2u_{\beta}}
$, which implies that $
\frac{U_{\lambda}^{2}}{w_{\beta}^{2}}\leqslant \frac{2L_{g,1}^{2}\lambda ^2}{w_{\beta}^{2}}\leqslant \frac{1}{u_{\beta}^{2}}
$. Then, we can upper bound the terms \(D_1\), \(D_2\) and \(D_3\) in \eqref{EQ-d0-d2-V} as follows:
\begin{equation}
D_1\leqslant p_1\lambda ^2, D_2\leqslant p_2\lambda ^2, D_3\leqslant p_3\lambda ^2,
\end{equation}
where   \(p_1\), \(p_2\) and \(p_3\) are given by \eqref{EQ-ppp}. Next,  by combining Lemma \ref{LE-consensus errors} and the inequality \eqref{EQ-d0-d2-V} and incorporating the potential function $V^k$  in
 \eqref{EQ-potential}, we further obtain the descent of the potential function $V^k$ as follows:
\begin{align}
&V^{k+1} \\
\leqslant& V^k-\alpha \left\| \nabla \Phi \left( \bar{x}^k \right) \right\| ^2 \nonumber\\
&+\mathcal{C}(\lambda, \alpha,\beta, C_{\rm ou}, C_{\rm in},\rho) \alpha +\mathcal{B}(\lambda,\alpha,\beta,\gamma,b_f, b_g,\rho) \alpha
\nonumber\\
&-\!C_1\alpha \left\| \bar{H}_{x}^{k} \right\| ^2\!-\!C_2\left\| \bar{y}^k-y^*\left( \bar{x}^k;\lambda \right) \right\| ^2\!-\!C_3\left\| \bar{z}^k-y^*\left( \bar{x}^k \right) \right\| ^2
\nonumber \\
&-\frac{C_4}{m}\left\| x^k-1_m\otimes \bar{x}^k \right\| ^2-\frac{C_5}{m}\left\| y^k-1_m\otimes \bar{y}^k \right\| ^2 \nonumber\\
&-\frac{C_6}{m}\left\| z^k-1_m\otimes \bar{z}^k \right\| ^2,  \nonumber
\end{align}
where $\mathcal{C}(\lambda, \alpha,\beta, C_{\rm ou}, C_{\rm in},\rho)$ and $\mathcal{B}(\lambda,\alpha,\beta,\gamma,b_f, b_g,\rho)$ are the penalty errors
and heterogeneity errors given by \eqref{EQ-p-error}  and  \eqref{EQ-b-error}, respectively, and
the terms $C_1$ to $C_6$ are given by:
\begin{align}
&\!\!\!C_1= \frac{1}{2}-\frac{24L_{y^*}^{2}L_{g,1}^{2}\lambda^2}{w_{\gamma}^{2}\gamma ^2}\alpha ^2-\frac{24U_{\lambda}^{2}L_{y*,\lambda}^{2}}{w_{\beta}^{2}\beta ^2}\alpha ^2,
\\
&\!\!\!C_2=12U_{\lambda}^{2}\alpha -\frac{108U_{\lambda}^{2}\alpha ^2d_3}{1-\rho}-\frac{24L_{g,1}^{2}\gamma ^2d_5}{1-\rho},
\\
&\!\!\!C_3=12L_{g,1}^{2}\lambda ^2\alpha -\frac{108U_{\lambda}^{2}\alpha ^2d_3}{1-\rho}-\frac{72U_{\lambda}^{2}\beta ^2d_4}{1-\rho},
\\
&\!\!\!C_4 \!= \! \frac{1-\rho}{4}d_3\!-\!\frac{108U_{\lambda}^{2}\alpha ^2d_3}{1-\rho} \!-\! \frac{72U_{\lambda}^{2}\beta ^2\!d_4}{1-\rho}\!-\!\frac{24L_{g,1}^{2}\gamma ^2d_5}{1-\rho},
\\
&\!\!\!C_5=\frac{1-\rho}{4}d_4-\frac{108U_{\lambda}^{2}\alpha ^2d_3}{1-\rho}-\frac{72U_{\lambda}^{2}\beta ^2d_4}{1-\rho},
\\
&\!\!\!C_6=\frac{1-\rho}{4}d_5-\frac{108L_{g,1}^{2}\alpha ^2\lambda ^2d_3}{1-\rho}-\frac{24L_{g,1}^{2}\gamma ^2d_5}{1-\rho}.
\end{align}
Note that,  the terms $C_1$ to $C_6$ are nonnegative, when the conditions \eqref{EQ-stepsize-alpha}-\eqref{EQ-stepsize-gamma} hold, which further gives
\begin{align}
V^{k+1}&\leqslant V^k-\alpha \left\| \nabla \Phi \left( \bar{x}^k \right) \right\| ^2   \label{EQ-Vk1-Vk}\\
&+\mathcal{C}(\lambda, \alpha,\beta, C_{\rm ou}, C_{\rm in},\rho)\alpha
+\mathcal{B}(\lambda,\alpha,\beta,\gamma,b_f, b_g,\rho) \alpha.  \nonumber
\end{align}
Summing the above inequality from \( k = 0 \) to \( K-1 \) yields the result in \eqref{EQ-the1}.
This completes the proof. {\hfill $\blacksquare$}

\colorbf{\textbf{2) Proof of Corollary  \ref{COR-rate}.}} To provide a tight analysis w.r.t. the condition number, we first need to analyze how the parameters  in \eqref{EQ-the1} depend on the condition number \(\kappa\). Specifically, by the definitions of $L$, $L_{y^*}$, $L_{y^*,\lambda}$, we have that $L=\mathcal{O}(\kappa^3)$, $L_{y^*}=\mathcal{O}(\kappa)$, $L_{y^*,\lambda}=\mathcal{O}(\kappa)$\footnote{{The symbol $\mathcal{O}(\cdot)$ used here hides constants that depend on the properties of the functions, emphasizing the dependency on  the condition number.}}. In addition, from \eqref{EQ-inner-errors}, \eqref{EQ-inner-penalty} and \eqref{EQ-potential}, we have that \(\frac{1}{w_{\gamma}} = \mathcal{O}(\kappa)\), \(\frac{1}{w_{\beta}} = \mathcal{O}\left(\frac{\kappa}{\lambda}\right)\), and  \(\frac{1}{u_{\beta}} = \mathcal{O}(\kappa)\), which implies that \(p_1 = \mathcal{O}(\kappa^2)\), \(p_2 = \mathcal{O}(\kappa^2)\), and \(p_3 = \mathcal{O}(\kappa^2)\) in \eqref{EQ-potential}. Furthermore, the coefficients \(d_1\) and \(d_2\) of the potential function \(V^0\) are derived as \(d_1 = \mathcal{O}\left(\kappa \lambda \frac{\alpha}{\beta}\right)\) and \(d_2 = \mathcal{O}\left(\kappa \lambda^2 \frac{\alpha}{\gamma}\right)\). When the variables \(\{x_i^0\}\), \(\{y_i^0\}\), and \(\{z_i^0\}\) are initialized as \(x_i^0 = x_j^0\), \(y_i^0 = y_j^0\), and \(z_i^0 = z_j^0\) for all \(i, j \in \mathcal{V}\), we additionally obtain \(V^0 = \mathcal{O}\left(\kappa \lambda \frac{\alpha}{\beta} \left\| \bar{y}^0 - y^*(\bar{x}^0; \lambda) \right\|^2 + \kappa \lambda^2 \frac{\alpha}{\gamma} \left\| \bar{z}^0 - y^*(\bar{x}^0) \right\|^2\right)\). To achieve the optimal \(\kappa\)-dependence, we set \(\alpha = \mathcal{O}(\kappa^{-3}K^{-\frac{2}{3}})\), \(\beta = \mathcal{O}(\kappa^{-1}K^{-\frac{1}{2}})\), \(\gamma = \mathcal{O}(K^{-\frac{1}{3}})\), and \(\lambda = \mathcal{O}(\kappa K^{\frac{1}{6}})\) such that the conditions \eqref{EQ-stepsize-alpha}-\eqref{EQ-stepsize-gamma} are satisfied and the term $V^0$ is sufficiently small. By combining  the fact that \(C_{\rm ou} = \mathcal{O}(\kappa^3)\) in \eqref{EQ-constant-symbols-cc} and substituting the chosen step sizes \(\alpha\), \(\beta\), \(\gamma\), and the penalty parameter \(\lambda\) into \eqref{EQ-the1}, we arrive at result \eqref{EQ-cor1-1}.
This completes the proof. {\hfill $\blacksquare$}

\colorbf{\textbf{3) Proof of Corollary  \ref{COR-cor2}.}}
 Our proof is split into three parts: \textbf{i)} First, if \(g_i(x, y) = -f_i(x, y)\) holds for all \(i \in \mathcal{V}\), it follows that \(y^*(x) = \arg \max _y (1/m) \sum_{i=1}^{m} f_i(x, y)\) in the problem \eqref{EQ-P1}. In this case, by the definition of the minimax problem, we can then conclude that the problem \eqref{EQ-P1} is equivalent to the problem \eqref{EQ-minmax} \cite{razaviyayn2020nonconvex}.

\textbf{ii)} Given the condition \(g_i(x, y) = -f_i(x, y)\) for all \(i \in \mathcal{V}\),
it holds that $y^*(x;\lambda)=\arg\min_y\max_z p(x,y,z;\lambda)=\arg\min_y\frac{1}{m}\sum_{i=1}^m (\left( 1-\lambda \right) f_i\left( x,y \right) +\lambda f_i\left( x,y^*\left( x \right) \right))$.
Additionally, by Assumption \ref{ASS-inner-level}(i), we know that  \(f_i(x, y)\) is strongly concave in \(y\). Therefore, when \(\lambda > 1\),  the function \(\frac{1}{m}\sum_{i=1}^m(1 - \lambda)f_i(x, y)\) is strongly convex in \(y\) and it holds that \(y^*(x, \lambda) = y^*(x)\) and \(p^*(x; \lambda) = \Phi(x)\).  This fact implies that, in the distributed minimax case, the problem \eqref{EQ-minmax} becomes equivalent to its penalty formulation $\min_{x,y}\max_z p(x,y,z;\lambda)$ and the inner-level penalty gap  reduces to zero, i.e.,
${\left\| y^*\left( x \right) -y^*\left( x;\lambda \right) \right\| ^2=0}$ and $\| \nabla \Phi \left( x \right) -\nabla p^*\left( x;\lambda \right) \|^2=0$ in Lemma \ref{LE-inner-penalty-error} with $C_{\rm in}=0$ and $C_{\rm ou}=0$. In addition, when $\lambda =1$, the penalty formulation exactly  reduces to the problem \eqref{EQ-minmax} and it holds that  ${\left\| y^*\left( x \right) -y^*\left( x;\lambda \right) \right\| ^2=0}$ and $\| \nabla \Phi \left( x \right) -\nabla p^*\left( x;\lambda \right) \|^2=0$ with $C_{\rm in}=0$ and $C_{\rm ou}=0$.

\textbf{iii)} When \(g_i(x, y) = -f_i(x, y)\) for all \(i \in \mathcal{V}\), and Assumptions \ref{ASS-outer-level}, \ref{ASS-network}, \ref{ASS-heterogeneity} and  Assumptions \ref{ASS-inner-level}(i), \ref{ASS-inner-level}(ii) hold, we can also derive Lemmas \ref{LE-descent}, \ref{LE-hypergradient-approx}, \ref{LE-inner-errors}, \ref{LE-inner-penalty}, \ref{LE-consensus errors} with the  Lipschitz constant of $\nabla \Phi(x)$ being \(L = L_{f,1}+L_{f,1}L_{y^*}\). Moreover,  if the conditions \eqref{EQ-stepsize-alpha}-\eqref{EQ-stepsize-gamma} are satisfied, the inequality \eqref{EQ-Vk1-Vk} holds with the penalty error $\mathcal{C}(\lambda, \alpha,\beta, C_{\rm ou}, C_{\rm in},\rho)$ being zero. Additionally, from the respective definitions of  \(w_{\gamma}\)  and \(w_{\beta}\) in \eqref{EQ-inner-errors} and \eqref{EQ-inner-penalty}, we can derive that the coefficients \(p_1\), \(p_2\), and \(p_3\) satisfy \(p_1 = \mathcal{O}(\kappa^2)\), \(p_2 = \mathcal{O}(\kappa^2)\), and \(p_3 = \mathcal{O}(\kappa^2)\) when \(\lambda\) is independent of the condition number and the number of iterations. Finally, by telescoping the inequality \eqref{EQ-Vk1-Vk} from 0 to \(K-1\), and noting that \(b_{f} = b_{g}\) in the term $\mathcal{B}(\lambda,\alpha,\beta,\gamma,b_f, b_g,\rho)$, we obtain:
\begin{equation}
\begin{aligned}\label{EQ-COR2-alpha}
& \frac{1}{K}\sum_{k=0}^{K-1}{\| \nabla \Phi \left( \bar{x}^k \right) \| ^2} \\
 \leqslant&
\mathcal{O}\big(\! \frac{V^0-V^{K-1}}{ K\alpha}\!\!+\! \!\frac{\kappa^2\alpha ^2b_{f}^{2}\!\!+\!\!\kappa^2\beta ^2b_{f}^{2}\!\!+\!\!\kappa^2\gamma ^2b_{f}^{2}}{\left( 1-\rho \right) ^2} \big).
\end{aligned}
\end{equation}
Furthermore, if the step sizes $\alpha$,  $\beta$, and $\gamma$ respectively are taken as $\alpha=\mathcal{O}(\kappa^{-2}K^{-{\frac{1}{3}}})$, $\beta=\mathcal{O}(K^{-{\frac{1}{3}}})$, $\gamma=\mathcal{O}(K^{-{\frac{1}{3}}})$ with $\lambda\alpha=\mathcal{O}(\kappa^{-2}K^{-{\frac{1}{3}}})$ and $\lambda\beta=\mathcal{O}(K^{-{\frac{1}{3}}})$,   we can conclude that the condition \eqref{EQ-stepsize-alpha}-\eqref{EQ-stepsize-gamma} are satisfied. Following a similar argument as in the proof of Corollary   \ref{COR-rate}, it holds that \(V^0 = \mathcal{O}(1)\) are independent
of the condition number and the number of iteration. Combining these results yields \eqref{EQ-COR2-2}. This completes the proof. {\hfill $\blacksquare$}

\section{Numerical Experiments}
In this section, we will provide two examples to verify the performance of the proposed algorithm.
\begin{figure}[ht]
\centering
\subfloat[]{
\begin{minipage}[ht]{0.5\linewidth}
\centering
\includegraphics[width=1\textwidth,height=0.15\textheight]{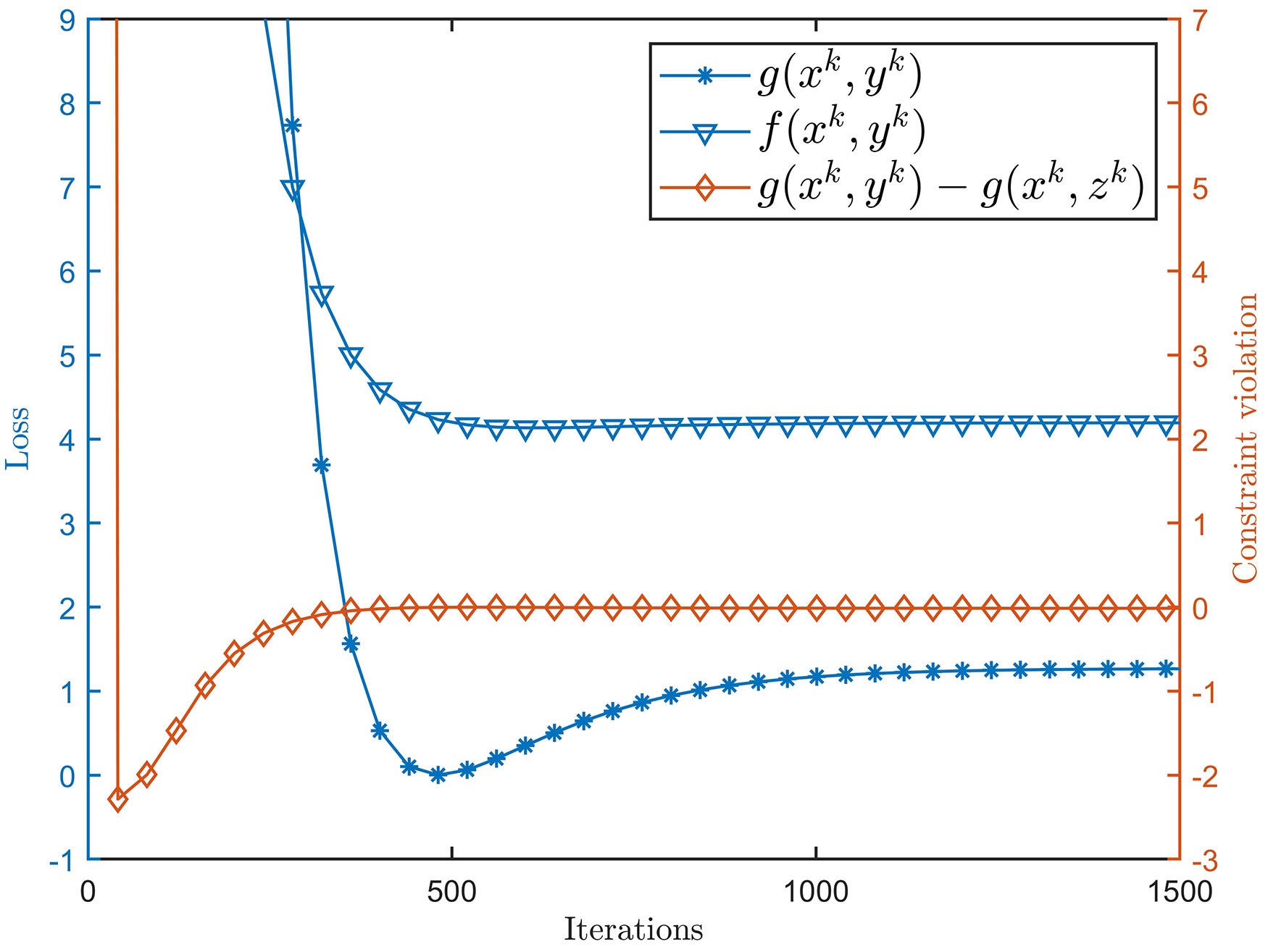} 
\end{minipage}%
}%
\centering
\subfloat[]{
\begin{minipage}[ht]{0.5\linewidth}
\centering
\includegraphics[width=1\textwidth,height=0.151\textheight]{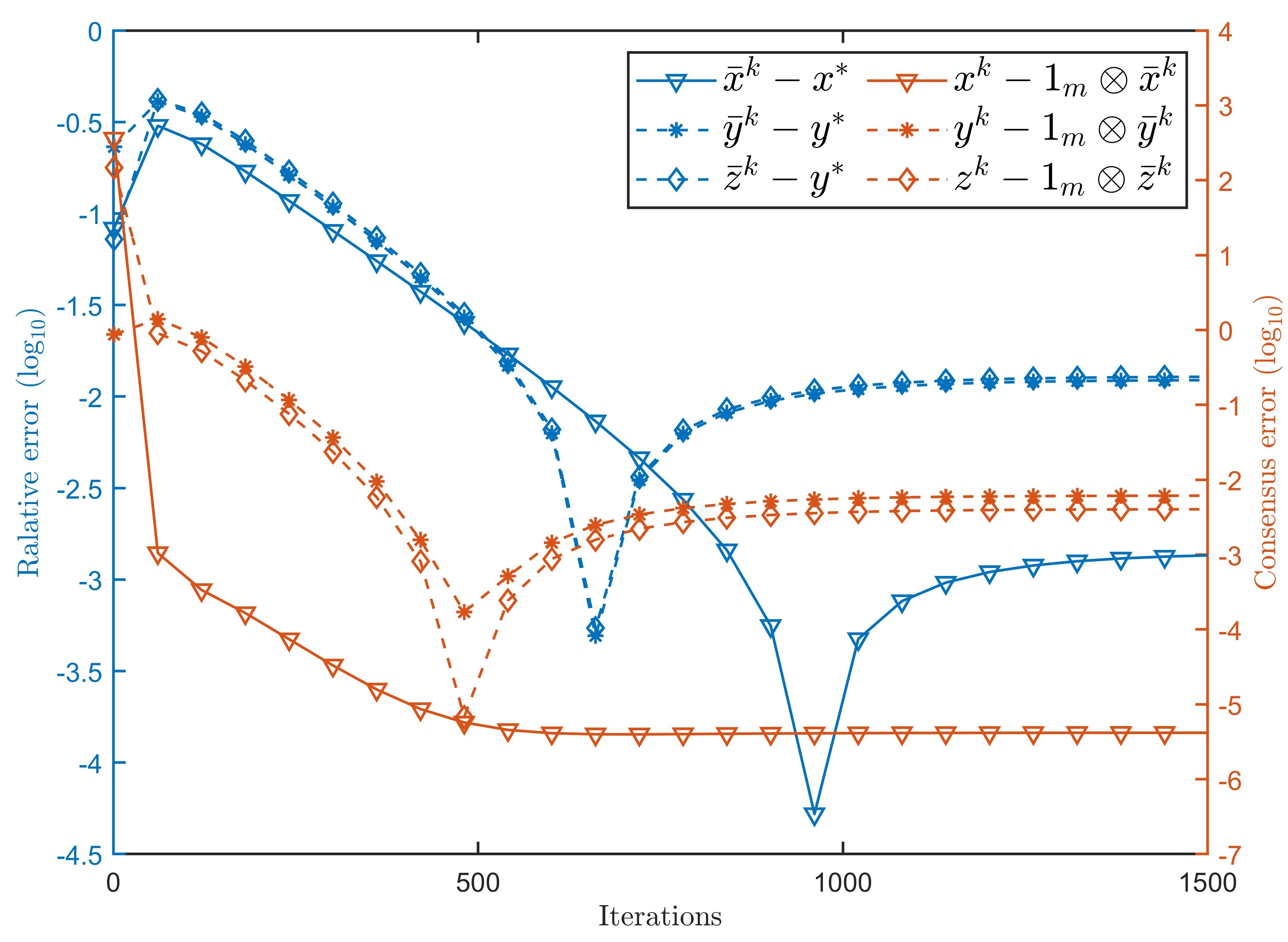}
\end{minipage}%
}%
\\
\subfloat[]{
\begin{minipage}[ht]{0.5\linewidth}
\centering
\includegraphics[width=1\textwidth,height=0.15\textheight]{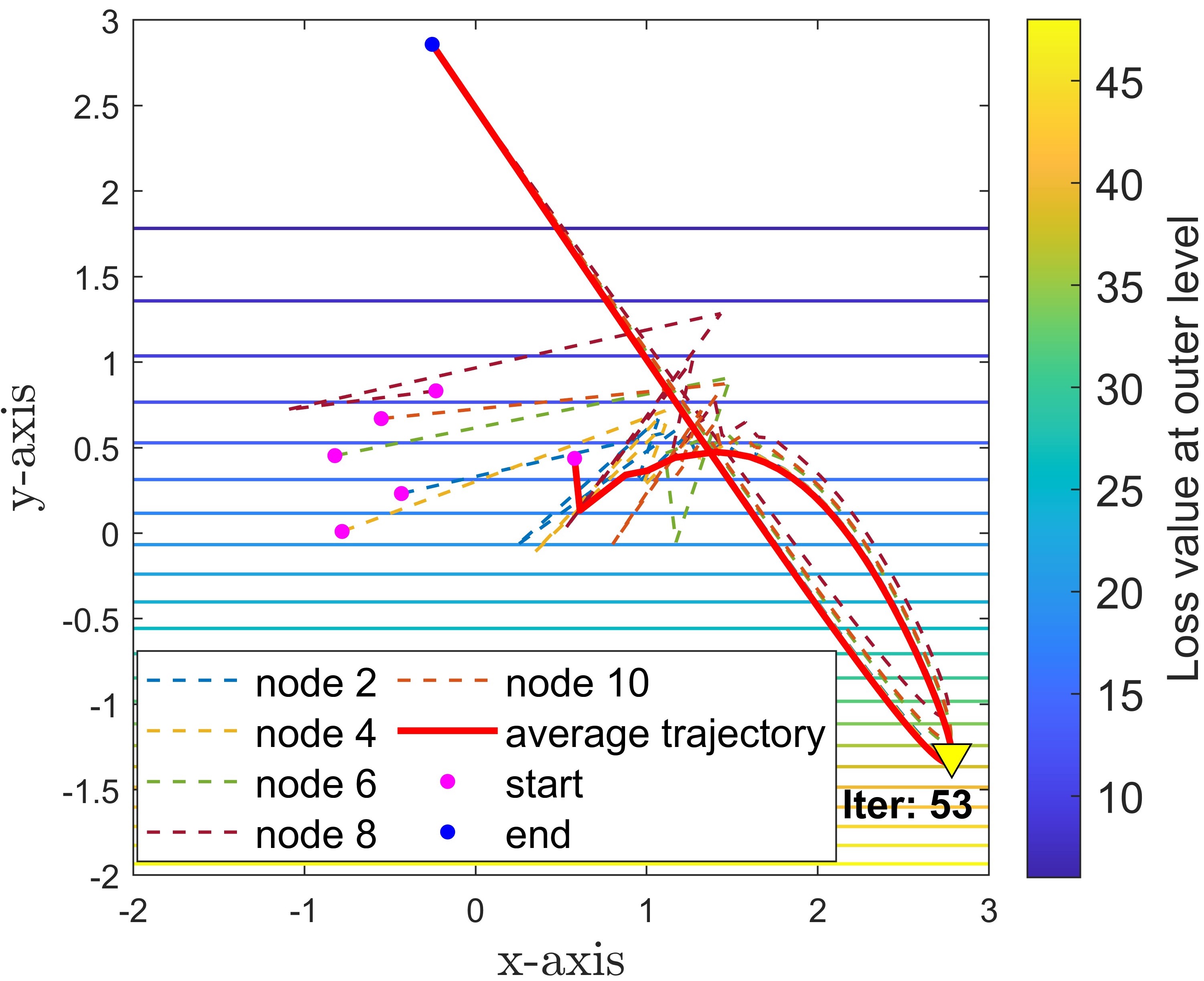} 
\end{minipage}%
}%
\centering
\subfloat[]{
\begin{minipage}[ht]{0.5\linewidth}
\centering
\includegraphics[width=1\textwidth,height=0.15\textheight]{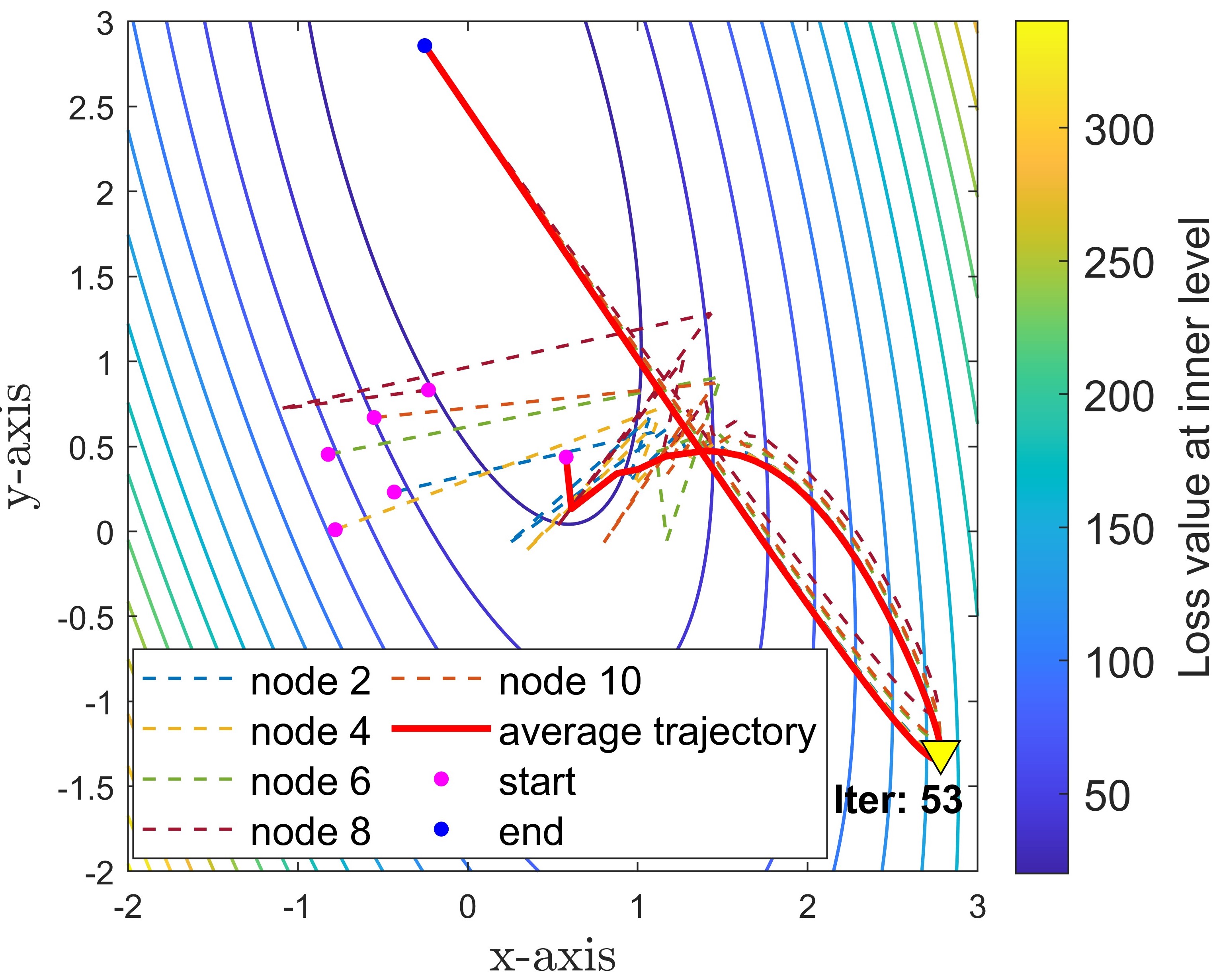}
\end{minipage}%
}%
\centering
\caption{The performance of the proposed {\ALGNAME} algorithm:  (a) The loss  \(f,g\) and constraint tolerance   with respect to iterations; (b) The optimal gap  and consensus errors  with respect to the number of iterations;  (c) and (d): The trajectories of the individual state $(x_i^k, y_i^k)$ and the average state $(\bar{x}^k, \bar{y}^k)$ within the contours of the outer- and inner-level objectives, respectively.}
\label{fig:case1}
\end{figure}
\subsection{Synthetic Example}
We begin by demonstrating the efficiency of the proposed algorithm using a toy example, where the local outer-level and inner-level functions are expressed as: $
f_i\left( x,y \right) =\frac{1}{2}\left( a_iy-b_i \right) ^2
$ and $
g_i\left( x,y \right) =\frac{1}{2}\left( c_ix+d_iy-e_i \right) ^2
$. As for data setting, we consider a network with $m=10$ nodes and a connectivity probability of $p_{\rm net}=0.7$, generated by the Erd\H{o}s-R\'enyi model. The parameters $\{a_i,b_i,c_i,d_i,e_i\}$ are set as follows: $a_i = 2$ for all $i = 1, \dots, m$; $c_i = 2$, $d_i = 2$, and $e_i = 10$ for $i = 1, \dots, 5$, and $c_i = 4$, $d_i = 4$, and $e_i = 10$ for $i = 6, \dots, m$; finally, $b_i = i$ for $i = 1, \dots, 10$. Under such parameters, the optimal solution to the problem \eqref{EQ-P1} can be calculated as $x^*=0.25$ and $ y^*=2.75$. The step sizes are set as $\alpha=0.0007$, $\beta=0.001$, $\gamma=0.01$ with $\lambda=20$. The performance of the {\ALGNAME} algorithm is presented in Fig.  \ref{fig:case1}, where the relative error and consensus error, labeled on the vertical axis,  are measured using the quantity $\log \|\cdot\|^2$.

As shown in Fig. \ref{fig:case1}a, both the losses  of the inner-level function \( g(x^k, y^k) = \frac{1}{m} \sum_{i=1}^m g_i(x_i^k, y_i^k) \)  and the outer-level function \( f(x^k, y^k) = \frac{1}{m} \sum_{i=1}^m f_i(x_i^k, y_i^k) \) decrease approximately at  the same rate. Additionally, the constraint \( g(x^k, y^k) - g(x^k, z^k) \leq 0 \) is always satisfied, and the convergence of the term $g(x^k, y^k) - g(x^k, z^k)$  is faster than that of the losses associated with both the inner-level and outer-level functions. These observations are consistent with the result in  Theorem \ref{COR-rate}  and Remark \ref{RE-revised-heter} under multiple-timescale step sizes.
On the other hand, the blue curves in Fig. \ref{fig:case1}b demonstrate that the {\ALGNAME} algorithm successfully converges to the optimal solution. The orange curves in Fig. \ref{fig:case1}b show that the consensus errors associated with the outer-level variables $\{x_i^k\}$ converge at a faster rate than those of the inner-level-related variables $\{y_i^k\}$ and $\{z_i^k\}$, which is consistent with our theoretical findings. Figs. \ref{fig:case1}c and \ref{fig:case1}d show the trajectories of \((x_i^k, y_i^k)\) and \((\bar{x}^k, \bar{y}^k)\) over iterations 0–1500. We observe that  with a random initialization, the trajectory \((x_i^k, y_i^k)\) asymptotically reaches consensus and quickly converges to the solution point once the consensus errors of the inner-level variables \(\{z_i^k\}\) and \(\{y_i^k\}\) begin to decrease around the 50th iteration, as shown in Fig. \ref{fig:case1}b. These observations further demonstrate the impact of the heterogeneity errors on the convergence and the crucial role of the inner-level heterogeneity.


\subsection{Hyperparameter Optimization}
In this example, we consider a distributed hyperparameter optimization problem in \( l_2 \)-regularized binary logistic regression for a binary classification scenario. Specifically, given the \( l_2 \)-regularization coefficient \(\eta \in \mathbb{R}^n\), a group of nodes aim to determine a logistic regression model \( y^*(x) \in \mathbb{R}^n \)  on  training datasets \(\{\mathcal{D}_i^{\rm train}\}\) at the inner level, while maximizing the performance of the model \( y^*(x) \) on validation datasets \(\{\mathcal{D}_i^{\rm val}\}\) by optimizing the coefficient \(\eta\), i.e.,
\begin{align}
&\min_{\eta \in \mathbb{R} ^n} \frac{1}{m}\sum_{i=1}^m{\sum\nolimits_{(s_{ij},b_{ij})\in \mathcal{D} _{i}^{\mathrm{val}}}{\varphi \left( y^*(\eta );b_{ij},s_{ij} \right)}} \label{EQ-hyper-opt}
\\
&\mathrm{s}.\mathrm{t}.\;  y^*(\eta )=\mathrm{arg}\min_{y\in \mathbb{R} ^n} \frac{1}{m}\sum_{i=1}^m{\sum\nolimits_{(s_{ij},b_{ij})\in \mathcal{D} _{i}^{\mathrm{train}}}{\varphi \left( y;b_{ij},s_{ij} \right)}} \nonumber \\
 &\quad\quad\quad\quad \quad\quad\quad\quad\quad\quad+y^{\mathrm{T}}\mathrm{diag}\{\mathrm{e}^{\eta}\}y,    \nonumber
\end{align}
where \( e^\eta = \operatorname{col} \{ e^{\eta_t} \}_{t=1}^n \), $
\varphi \left( y;b,s \right) =\log\mathrm{(}1+e^{-(bs^{\mathrm{T}}y)})
$, and \( (s_{ij}, b_{ij}) \) denotes the \( j \)-th sample at node \( i \), with \( s_{ij} \in \mathbb{R}^n \) being the feature vector and \( b_{ij} \in \{1,-1\} \) the associated label.

We conduct experiments on a connected network of $m=10$ nodes with the MNIST dataset, using 4000  samples (digits `1' and `3', with each having 784 features) for training. First of all, we examine our {\ALGNAME} algorithm in following three cases, where the testing accuracy shown in all figures is computed using the average of the local variables $\{y_i\}$.

\begin{figure}[ht]
\centering
\begin{minipage}[ht]{0.7\linewidth}
\centering
\includegraphics[width=0.95\textwidth,height=0.20\textheight]{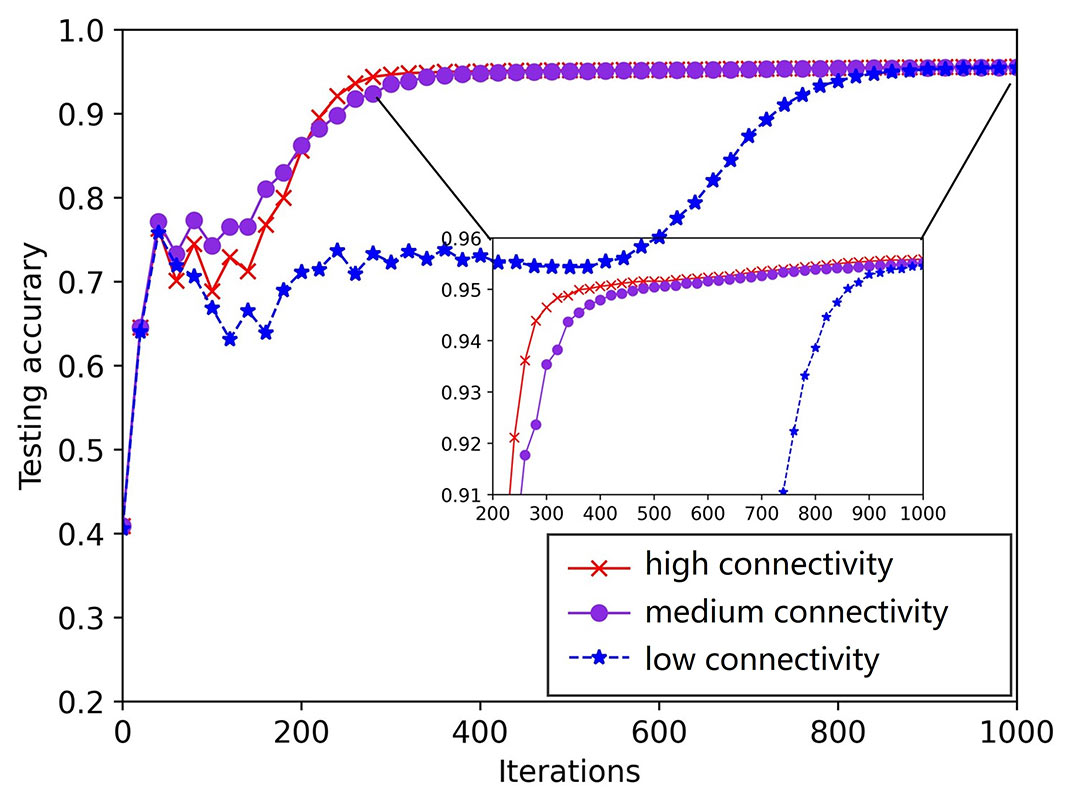} 
\end{minipage}%
\centering
\caption{The performance of the proposed {\ALGNAME} algorithm under networks with different connectivity levels, where high, medium, and low connectivity correspond to spectral gaps of \(\rho = 0.274\), \(\rho = 0.644\), and \(\rho = 0.923\), respectively.}
\label{fig:case2-heter}
\end{figure}

 \colorbf{\textbf{i)} \textbf{Different network connectivity.}} We evaluate our algorithm on three types of networks with spectral gaps of \(\rho = 0.274\), \(\rho = 0.644\), and \(\rho = 0.923\). The training and validation sets are randomly assigned, ensuring an equal distribution of samples across nodes. The experimental results are shown in Fig. \ref{fig:case2-heter}, where we observe that as network connectivity decreases, the number of iterations required to achieve a certain testing accuracy increases.
Notably, when the network transitions from a relatively well-connected state to a weakly connected one, this increase becomes more pronounced. These results verify our theoretical findings that network connectivity impacts algorithm convergence,  particularly during the initial stages of the process.

 \colorbf{\textbf{ii)} \textbf{Different step sizes and penalty parameters}}. We examine the impact of step sizes and the penalty parameter on the algorithm's performance. The results, presented in Table \ref{TAB-stepsize}, show testing accuracy and consensus errors under eight typical selections. In particular, as shown in the columns 1 and 2, when \(\gamma\) is small, a small penalty \(\lambda\) guarantees convergence but with low accuracy, while a larger \(\lambda\) leads to non-convergence. The columns 2-4 shows that, under a large penalty \(\lambda\), increasing \(\gamma\) ensures convergence with higher accuracy but results in larger consensus errors. A selection with the time-scale separation as \(\gamma > \beta > \alpha\) guarantees convergence with higher accuracy, as indicated in columns 4-7 of  Table \ref{TAB-stepsize}. A larger value of  \(\lambda\) improves accuracy but results in larger consensus errors (c.f., columns 5 and 8 in Table~\ref{TAB-stepsize}). These results demonstrate the trade-off between the selection of step sizes and the penalty parameter, as well as the convergence guarantee under time-scale separation, validating the theoretical findings in Theorem \ref{TH-the1} and Corollary \ref{COR-rate}. In practice, for selecting the step sizes and penalty parameter, one can first set an expected value for the products $\alpha\lambda$ and $\beta\lambda$, along with a near-equal order of the value for $\gamma$. Then, $\alpha$ (or $\beta$) and $\lambda$ can be determined based on the value of the product $\alpha\lambda$ (or $\beta\lambda$), according to the observed trade-off.

\begin{table}[!ht]
\centering
\footnotesize{
\caption{The performance of the proposed {\ALGNAME} algorithm under different step sizes and penalty parameters.}
\label{TAB-stepsize}
\setlength{\tabcolsep}{0.070cm}{}
\begin{threeparttable}
\begin{tabular}{ccccccccc} \hline
   No.  &1   &2  &3  &4 &5 &6 &7 &8\\  \hline
   $\alpha$  &\scriptsize{0.0001}  &\scriptsize{0.0001}  &\scriptsize{0.0001} &\scriptsize{0.0001}   &\scriptsize{0.0001} &\scriptsize{0.001}    &\scriptsize{0.0001}  &\scriptsize{0.0001}\\
   $\beta$  &\scriptsize{0.0005}  &\scriptsize{0.0005} &\scriptsize{0.0005} &\scriptsize{0.0005}   &\scriptsize{0.001}  &\scriptsize{0.0005} &\scriptsize{0.001} &\scriptsize{0.001}\\
   $\gamma$  &\scriptsize{0.001}    &\scriptsize{0.001}    &\scriptsize{0.005}   &\scriptsize{0.02}   &\scriptsize{0.02}  &\scriptsize{0.01} &\scriptsize{0.001} &\scriptsize{0.02}\\
   $\lambda$  &\scriptsize{10}  &\scriptsize{100}  &\scriptsize{100}  &\scriptsize{100} &\scriptsize{100}  &\scriptsize{100} &\scriptsize{100} &\scriptsize{50} \\
   Acc.  &0.8896   &\scriptsize{N.A.}  &\scriptsize{0.9550} &\scriptsize{0.9552} &\scriptsize{0.9613}  &\scriptsize{N.A.} &\scriptsize{N.A.} &\scriptsize{0.9555} \\
   ${\rm CE}_x$  &\scriptsize{9.428E-09}  &\scriptsize{N.A.}   &\scriptsize{1.527E-06} &\scriptsize{1.596E-06} &\scriptsize{1.780E-05} &\scriptsize{N.A.} &\scriptsize{N.A.} &\scriptsize{3.896E-07}\\
   ${\rm CE}_y$  &\scriptsize{3.627E-03}  &\scriptsize{N.A.}   &\scriptsize{4.598E-02} &\scriptsize{4.815E-02} &\scriptsize{4.127E-01} &\scriptsize{N.A.} &\scriptsize{N.A.} &\scriptsize{4.920E-02}\\
   ${\rm CE}_z$  &\scriptsize{4.296E-04}  &\scriptsize{N.A.}  &\scriptsize{1.842E-03} &\scriptsize{1.214E-02} &\scriptsize{4.356E-03} &\scriptsize{N.A.} &\scriptsize{N.A.} &\scriptsize{1.216E-02}\\  \hline
\end{tabular}
\footnotesize{*Acc. refers to the testing accuracy after running 800 iterations; ${\rm CE}_x$, ${\rm CE}_y$, and ${\rm CE}_z$ refer to the consensus errors with respect to $\{x_i^k\}$, $\{y_i^k\}$, and $\{z_i^k\}$, respectively; N.A. indicates that the result is not available due to non-convergence.}
\end{threeparttable}
}
\end{table}

\begin{figure}[ht]
\centering
\subfloat[]{
\begin{minipage}[ht]{0.5\linewidth}
\centering
\includegraphics[width=0.95\textwidth,height=0.15\textheight]{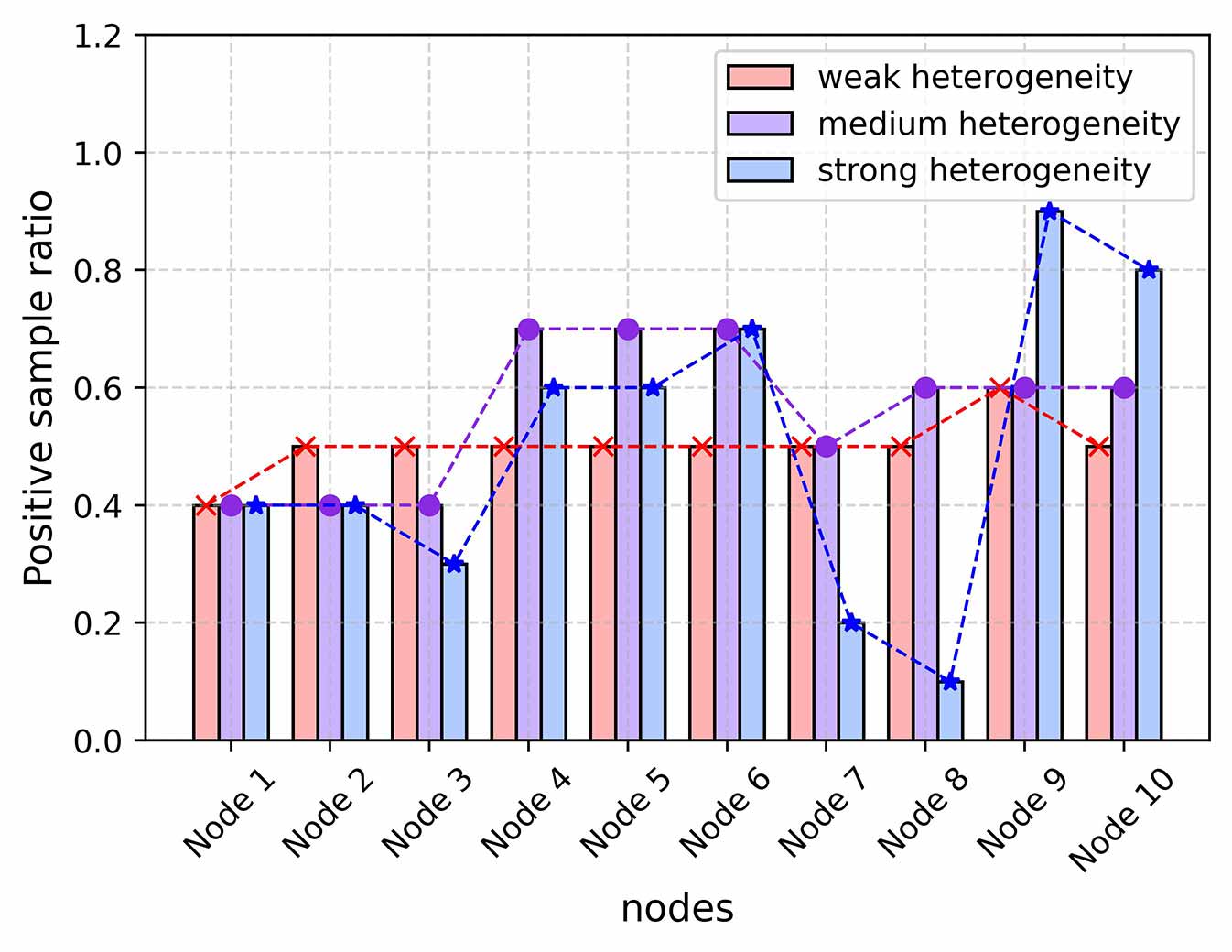} 
\label{fig:sample distribution}
\end{minipage}%
}%
\centering
\subfloat[]{
\begin{minipage}[ht]{0.5\linewidth}
\centering
\includegraphics[width=0.95\textwidth,height=0.15\textheight]{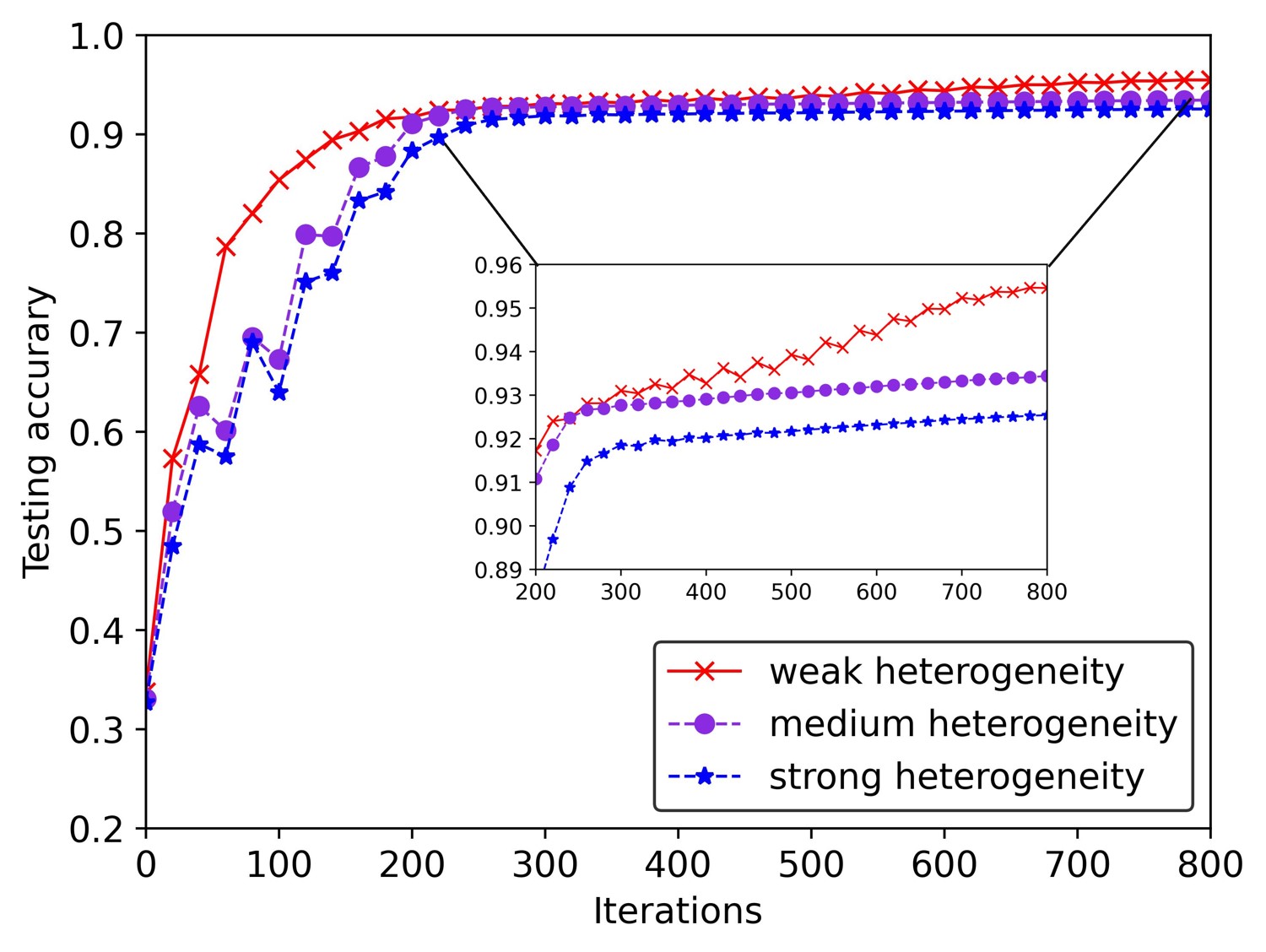}
\label{fig:heter-result}
\end{minipage}%
}%
\centering
\caption{The performance of the proposed AHEAD algorithm under different node heterogeneity: (a) Data distribution; (b) Testing accuracy.}
\label{fig:case2-different}
\end{figure}

\colorbf{\textbf{
\textbf{iii)}  Different node heterogeneity.}}  We test the algorithm's performance under three types of node heterogeneity. In this case, for each type, an equal number of samples is assigned to each node, with different positive sample ratios for each distribution. The detailed sample distribution for each type of heterogeneity is provided in Fig. \ref{fig:sample distribution}, and the experiment result is presented in  Fig. \ref{fig:heter-result}. It can be seen from Fig. \ref{fig:heter-result} that as the degree  of node heterogeneity decreases, both convergence speed and accuracy improve, which aligns with our theoretical results regarding the impact of the heterogeneity.

\begin{figure}[ht]
\centering
\subfloat[]{
\begin{minipage}[ht]{0.5\linewidth}
\centering
\includegraphics[width=0.95\textwidth,height=0.15\textheight]{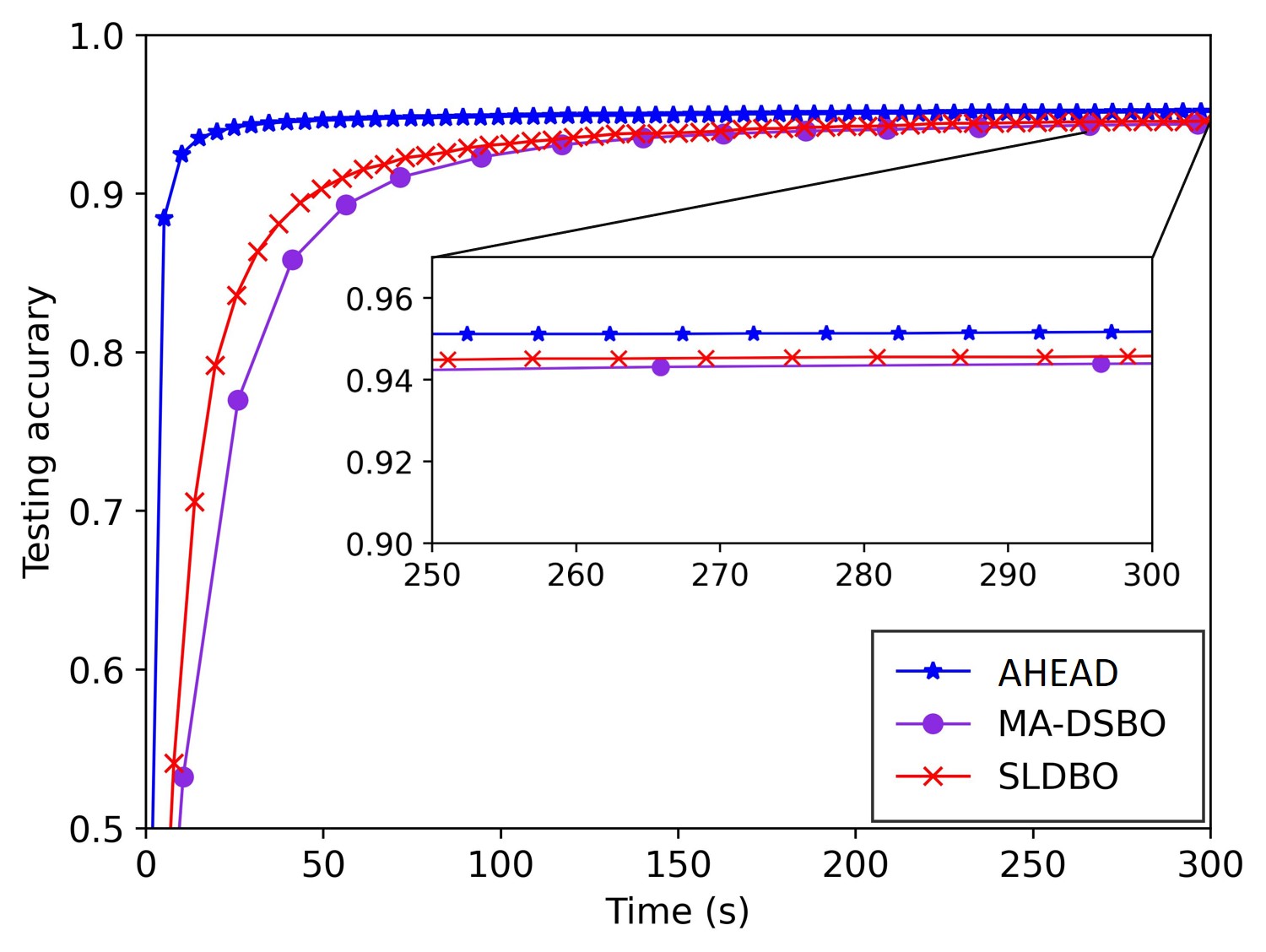} 
\end{minipage}%
}%
\centering
\subfloat[]{
\begin{minipage}[ht]{0.5\linewidth}
\centering
\includegraphics[width=0.95\textwidth,height=0.15\textheight]{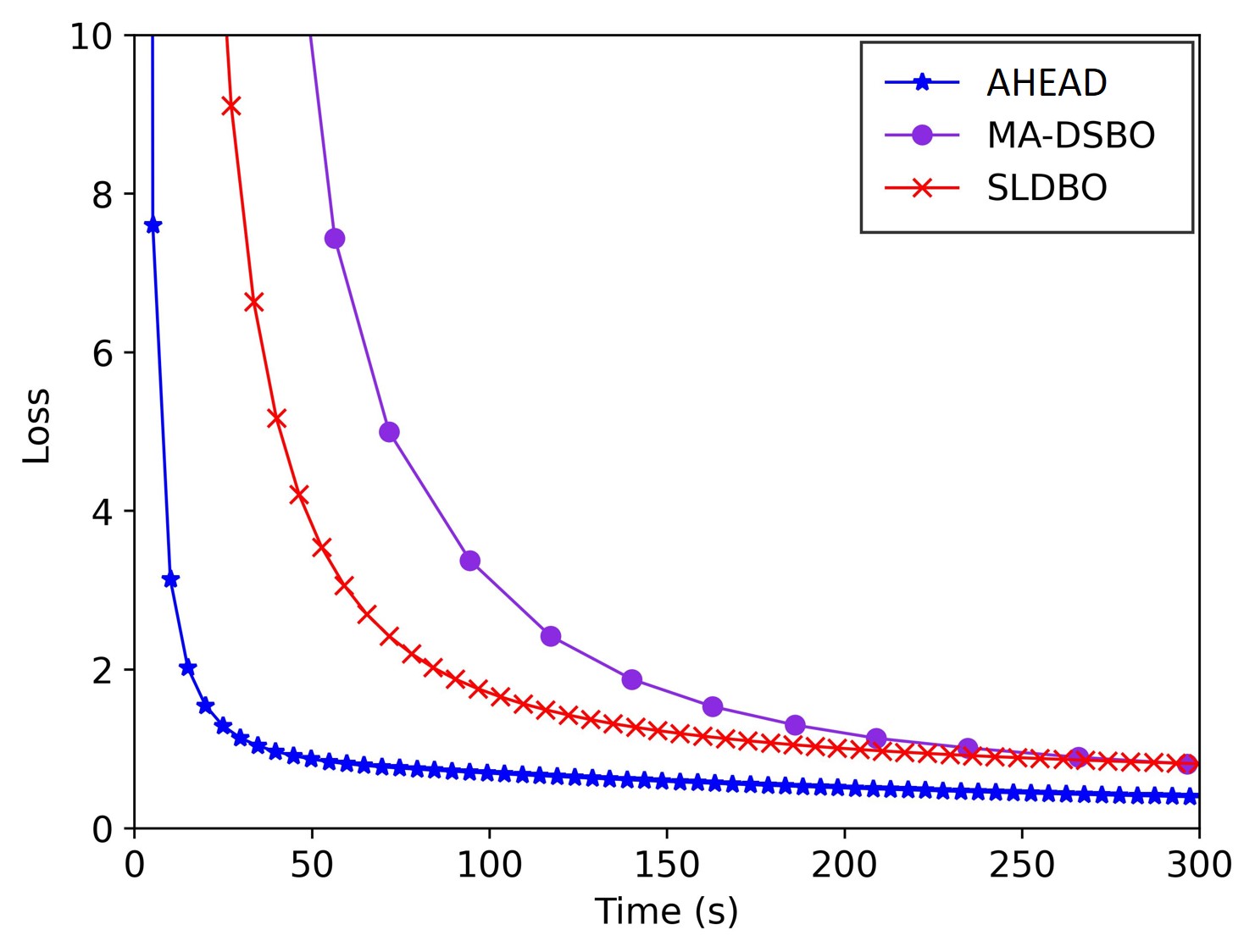}
\end{minipage}%
}%
\centering
\caption{Performance comparison of  MA-DSBO, SLDBO, and the proposed {\ALGNAME} algorithm: (a) Testing  accuracy; (b) Training loss.}
\label{fig:case2}
\end{figure}

To further verify the effectiveness of the proposed algorithm, we compare it with the representative Hessian-based distributed methods, including SLDBO \cite{dong2023single} and MA-DSBO \cite{chen2023decentralized} algorithm, in terms of computational time.  Specifically, SLDBO is a loopless Hessian-based algorithm that employs GT, while MA-DSBO is a double-loop Hessian-based algorithm that adds computational loops to estimate the inner-level solution and the Hessian inverse more accurately. To ensure a fair comparison, we adapt MA-DSBO by replacing minibatch gradients with full gradients, so all algorithms operate in the same deterministic setting. The training and validation sets are randomly assigned in each node with an equal number of samples.
The step sizes for both our algorithm and the baseline algorithms are manually tuned to achieve optimal performance. The experimental results are presented in Fig. \ref{fig:case2}. As shown in Fig. \ref{fig:case2}, the proposed {\ALGNAME} algorithm exhibits a significant advantage in computational time.

\section{Conclusions}
In this paper, we have proposed an approximated
 minimax formulation for the DBO problem with a coupled inner-level subproblem. Moreover, we have proposed a new distributed algorithm (AHEAD)  with multiple-timescale updates and  loopless structures.  The  algorithm features a  Hessian-free property, enjoying significant computational advantages in the large-scale problem where the solution accuracy is not essential.
We  have  proved that the algorithm can achieve a sublinear convergence for nonconvex-strongly-convex cases, explicitly characterizing the impact of  the penalty parameter, node heterogeneity, network connectivity, and  condition number on the convergence performance, relying on a weaker assumption on first-order heterogeneity.
Moreover, we have provided a  detailed  analysis for the distributed minimax problem as a special case. Two numerical experiments were conducted to demonstrate the superiority of the algorithm and validate the theoretical results. In our future work,  a dynamic  update scheme can be incorporated for the penalty parameter to further enhance the algorithm's flexibility. Extending the algorithm to handle more general inner-level functions, possibly with nonconvexity, nonsmoothness, stochasticity and constraints, would also be an interesting direction.

\appendix
\subsection{Proof of Lemma \ref{LE-descent}} \label{S-p1}
By the recursion \eqref{EQ-alg-c} and the definition of $H_x^k$ in \eqref{EQ-HHH}, it follows that $
\bar{x}^{k+1}=\bar{x}^k-\alpha \bar{H}_{x}^{k}$. Then, according to the smoothness of $\Phi$ in Proposition \ref{PRO-PHI}, one can obtain
\begin{align}
&\Phi \left( \bar{x}^{k+1} \right) \nonumber
\\
\leqslant &
 \Phi \left( \bar{x}^k \right) + \left< \nabla \Phi \left( \bar{x}^k \right) ,x^{k+1}-x^{k} \right> +\frac{L}{2}\left\| x^{k+1}-x^{k}\right\| ^2  \nonumber\\
\leqslant &
 \Phi \left( \bar{x}^k \right) -\alpha \left< \nabla \Phi \left( \bar{x}^k \right) ,\bar{H}_{x}^{k} \right> +\frac{L}{2}\alpha ^2\left\| \bar{H}_{x}^{k} \right\| ^2
\\
\leqslant& \Phi \left( \bar{x}^k \right) \!+\!\frac{\alpha}{2}{ \left\| \nabla \Phi \!\left( \bar{x}^k \right) \!-\!\bar{H}_{x}^{k} \right\| ^2}-\frac{\alpha}{2}\left\| \nabla \Phi \left( \bar{x}^k \right) \right\| ^2\!-\!\frac{\alpha}{4}\left\| \bar{H}_{x}^{k} \right\| ^2,   \nonumber
\end{align}
where the last step  utilizes the identity $-2a^{\rm T}b=\|a-b\|^2-\|a\|^2-\|b\|^2$ for any vectors $a$ and $b$, and subsequently applies the condition that $\alpha\leqslant\frac{1}{2L}$. This completes the proof. {\hfill $\blacksquare$}

\subsection{Proof of Lemma \ref{LE-inner-penalty-error}} \label{S-p1-0}
By the definition of $f(x,y)=\frac{1}{m}\sum_{i=1}^mf_i(x,y)$ and  $g(x,y)=\frac{1}{m}\sum_{i=1}^mg_i(x,y)$, the expression of $p(x, y, z; \lambda)$ can be compactly rewritten as \( p(x, y, z; \lambda) = f(x,y)+\lambda(g(x,y)-g(x,z))\).  Leveraging the gradient Lipschitz continuity of \( f(x, y) \) and the strong convexity of \( g(x, y) \) in \( y \), as imposed by Assumptions \ref{ASS-outer-level} and \ref{ASS-inner-level}, it follows that for \( \lambda > \frac{2 L_{f,1}}{\mu_g} \), \( p(x, y, z; \lambda) \) is \( \frac{\lambda \mu_g}{2} \)-strongly convex in \( y \). Moreover, since \( p(x, y, z; \lambda) \) is strongly concave in \( z \) and \( y^*(x) = \arg\max_z p(x, y, z; \lambda) \), it holds that  \( \min_y \max_z p(x, y, z; \lambda)=\min_y p(x, y, y^*(x); \lambda) \), and the solution \( y^*(x; \lambda) \) to the problem \( \min_y \max_z p(x, y, z; \lambda) \) is unique. Then, utilizing the strong convexity of $g(x,y)$ in $y$, we have the following inequality for the inner-level penalty gap:
\begin{align}
&\mu_g\|y^*(x)-y^*(x;\lambda)\| \label{EQ-ystar-ylambda-bound}\\
= & \|\nabla_yg(x,y^*(x;\lambda))\|=\|\frac{\nabla_y f(x,y^*(x;\lambda))}{\lambda}\| , \nonumber \\
\leqslant&  \frac{ \|\nabla_yf(x,y^*(x))-\nabla_y f(x,y^*(x;\lambda))\|}{\lambda} + \frac{\|\nabla_yf(x,y^*(x))\|}{\lambda} \nonumber\\
\leqslant & \frac{L_{f,1}}{\lambda}\|y^*(x)-y^*(x;\lambda)\|+\frac{C_{f,y}}{\lambda},
\nonumber
\end{align}
where the second step follows from the fact that \( \nabla_y p(x, y^*(x;\lambda),y^*(x); \lambda) = \nabla_y f(x, y^*(x;\lambda)) + \lambda \nabla_y g(x, y^*(x;\lambda)) = 0 \); the final step uses the Lipschitz continuity of $\nabla f$ and the boundness of \( \nabla_y f(x, y^*(x)) \), as imposed by Assumptions \ref{ASS-outer-level}. By combining the condition \(\lambda > \frac{2L_{f,1}}{\mu_g}\), we obtain \(\frac{L_{f,1}}{\lambda} < \frac{\mu_g}{2}\). Substituting this into \eqref{EQ-ystar-ylambda-bound} yields \(\|y^*(x) - y^*(x;\lambda)\| \leqslant \frac{2C_{f,y}}{\mu_g \lambda}=\frac{C_{\rm in}}{\lambda}\).

Next, we analyze the outer-level penalty gap. Recall that $p^*\left( x;\lambda \right)= p(x,y^*(x;\lambda),y^*(x);\lambda)$. By combining the optimality conditions for  $y^*(x;\lambda)$ and $y^*(x)$, we obtain $\nabla_y f(x, y^*(x;\lambda)) + \lambda \nabla_y g(x, y^*(x;\lambda)) = 0, \quad \text{and} \quad \nabla_y g(x, y^*(x)) = 0$, which further implies that $\nabla p^*(x;\lambda)=\nabla_xf(x,y^*(x;\lambda))+\lambda (\nabla_xg(x,y^*(x;\lambda))-\nabla_x g(x,y^*(x)))$.
By the expression for the hypergradient $\nabla \Phi(x)$  in  \eqref{EQ-hypergradient}, together with the expression of $\nabla p^*(x;\lambda)$,
the difference $\nabla \Phi \left( x \right) -\nabla p^*\left( x;\lambda \right)$ can be written as:
\begin{align}
&\nabla \Phi \left( x \right) -\nabla p^*\left( x;\lambda \right) \label{EQ-phi-p-expression} \\
= &\nabla _xf\left( x,y^*\left( x \right) \right) -\nabla _xf\left( x,y^*\left( x;\lambda \right) \right) \nonumber
\\
&+\nabla y^*\left( x \right) \left( \nabla _yf\left( x,y^*\left( x \right) \right) -\nabla _yf\left( x,y^*\left( x;\lambda \right) \right) \right)  \nonumber
\\
&-\lambda \nabla y^*\left( x \right) \pi _1-\lambda \pi _2,    \nonumber
\end{align}
where $\pi _1=\nabla _yg( x,y^*( x;\lambda ) ) -\nabla _yg( x,y^*( x ) ) -\nabla _{yy}^{2}g( x,y^*( x ) ) ( y^*( x;\lambda )\! -\!y^*( x ) )$ and $\pi _2=\nabla _xg( x,y^*( x;\lambda ) )$ $-\nabla _xg( x,y^*( x ) )$ $ -\nabla _{xy}^{2}g( x,y^*( x ) ) ( y^*( x;\lambda ) -y^*( x ) )$. By the Lipsthiz continuity of $\nabla_{xy}^2 g$ and $\nabla_{yy}^2 g$ imposed in Assumption \ref{ASS-inner-level}, we have that $
\| \pi _1 \| \leqslant \frac{L_{g,2}}{2}\| y^*( x ) -y^*( x;\lambda ) \|$ and $ \| \pi _2 \| \leqslant \frac{L_{g,2}}{2}\| y^*( x ) -y^*( x;\lambda ) \|
$. In addition, we have  $\| \nabla y^*( x ) \|=\|\nabla_{xy}^2(x,y^*(x))[\nabla_{yy}^2(x,y^*(x))]^{-1}\| \leqslant \frac{L_{g,1}}{\mu _g}$. Next, by applying the triangle inequality and the Cauchy-Schwarz inequality, along with the Lipschitz continuity of $\nabla_x f$ and $\nabla_x g$, it follows from \eqref{EQ-phi-p-expression} that
\begin{align}
&\| \nabla \Phi ( x ) -\nabla _xp( x,y;\lambda ) \| \label{EQ-phi-p-norm}\\
\leqslant &( 1\!+\!\frac{L_{g,1}}{\mu _g} ) \!\| y^*( x ) \!-\!y^*( x;\!\lambda ) \!\| ( L_{f,1}\!+\!\lambda \frac{L_{g,2}}{2}\| y^*( x ) \!-\!y^*( x;\!\lambda ) \| ). \nonumber
\end{align}
Then, by knowing the fact that \(\|y^*(x) - y^*(x;\lambda)\| \leqslant \frac{C_{\rm in}}{\lambda}\), we obtain the desired result. This completes the proof. {\hfill $\blacksquare$}

\subsection{Proof of Lemma \ref{LE-hypergradient-approx}} \label{S-p2}
Recall that  $p^*(x;\lambda)=p(x,y^*(x;\lambda),y^*(x);\lambda)$  in Lemma \ref{LE-inner-penalty-error}. Then, by employing the outer-level penalty gap in \eqref{LE-inner-penalty-error}, we can bound the term ${\left\| \nabla \Phi \left( \bar{x}^k \right) -\bar{H}_{x}^{k} \right\| ^2}$ as follows:
\begin{align}
& {\left\| \nabla \Phi \left( \bar{x}^k \right) -\bar{H}_{x}^{k} \right\| ^2}\label{EQ-phi-penalty}
\\
\leqslant & 2\left\| \nabla \Phi \left( \bar{x}^k \right) -\nabla p^*\left( \bar{x}^k;\lambda  \right) \right\| ^2+2\left\| \nabla p^*\left( \bar{x}^k;\lambda  \right) -\bar{H}_{x}^{k} \right\| ^2   \nonumber.
\end{align}
Note that the  term $\| \nabla p^*\left( \bar{x}^k;\lambda  \right) -\bar{H}_{x}^{k} \| ^2$ can be bounded as follows:
\begin{align}
 &\left\| \nabla p^*\left( \bar{x}^k;\lambda \right) -\bar{H}_{x}^{k} \right\| ^2 \label{EQ-penalty-star-h}
\\
\leqslant& 3\left\| J_n\nabla _xF\left( 1_m\otimes \bar{x}^k,1_m\!\otimes\! y^*\left( \bar{x}^k;\lambda  \right) \!\right) -J_n\nabla _xF\left( x^k,y^k \right) \right\| ^2 \nonumber
\\
&\!+3\!\lambda ^2\!\left\|\! J_n\!\nabla _xG\!\left( 1_m\!\otimes\! \bar{x}^k,1_m\!\otimes\! y^*\!\left( \bar{x}^k;\lambda  \right) \!\right) \!-\!J_n\!\nabla _xG\!\left( x^k,y^k \right) \right\| ^2 \nonumber
\\
&\!+3\lambda ^2\!\left\|\! J_n\nabla _xG\!\left( 1_m\otimes \bar{x}^k,1_m\otimes y^*\!\left( \bar{x}^k \right) \!\right)\! -\!J_n\!\nabla _xG\left( x^k,z^k \right) \right\| ^2 \nonumber
\\
\leqslant&6U_{\lambda}^2 \| \bar{y}^k\!-\!y^*( \bar{x}^k;\lambda  ) \| ^2\!+\!6L_{g,1}^{2}\lambda ^2\| \bar{z}^k\!-\!y^*( \bar{x}^k ) \| ^2 \nonumber
\\
&+6U_{\lambda}^2 \frac{1}{m}\left\| x^k-1_m\otimes \bar{x}^k \right\| ^2
+6U_{\lambda}^2 \frac{1}{m}\left\| y^k-1_m\otimes \bar{y}^k \right\| ^2 \nonumber \\
&+6L_{g,1}^{2}\lambda ^2\frac{1}{m}\left\| z^k-1_m\otimes \bar{z}^k \right\| ^2,   \nonumber
\end{align}
where the first inequality follows from  the definition of $H_x^k$ in \eqref{EQ-HHH} and the triangle inequality, and the second inequality uses the Lipschitz continuity of $\nabla_x g_i$ and $\nabla_x f_i$. Then, substituting  the inequalities \eqref{EQ-approx-error} and  \eqref{EQ-penalty-star-h} into the
 inequality \eqref{EQ-phi-penalty} yields the desired result \eqref{EQ-phi-hx-bar}.
 This completes the proof. {\hfill $\blacksquare$}

\subsection{Proof of Lemma \ref{LE-inner-errors}} \label{S-p3}
Using the term $y^*(\bar x^k)$ and employing Young's inequality,
$\left\| \bar{z}^{k+1}-y^*\left( \bar{x}^{k+1} \right) \right\| ^2$ can be bounded as follows:
\begin{align}
\!\left\| \bar{z}^{k+1}\!-\!y^*\left( \bar{x}^{k+1} \right) \right\| ^2
\leqslant &\left( 1+\gamma w_{\gamma} \right) \underset{\triangleq A_{1}^{z}}{\underbrace{{ \left\| \bar{z}^{k+1}-y^*\left( \bar{x}^k \right) \right\| ^2}}} \label{EQ-zbar-ystar}  \\
&+( 1+\frac{1}{\gamma w_{\gamma}} ) \underset{\triangleq A_{2}^{z}}{\underbrace{{\left\| y^*\left( \bar{x}^k \right) -y^*\left( \bar{x}^{k+1} \right) \right\| ^2}}}.   \nonumber
\end{align}
Next, we  combine the recursion \eqref{EQ-alg-a} with the definition of $H_z^k$ in \eqref{EQ-HHH}, and further bound the term $A_1^z$ in \eqref{EQ-zbar-ystar}:

\begin{align}
A_{1}^{z}=&\left\| \bar{z}^k-\gamma \bar{H}_{z}^{k}-y^*\left( \bar{x}^k \right) \right\| ^2 \label{EQ-A1Z}
\\
\leqslant &  \left( 1+\gamma w_{\gamma} \right) { \!\left\| \bar{z}^k \!- \!\gamma J\nabla _yG( 1_m\!\otimes\! \bar{x}^k,1_m\!\otimes\! \bar{z}^k ) -y^*\left( \bar{x}^k \right) \right\| ^2} \nonumber
\\ &\!\!\!\!\!\! +(\!1+\!\frac{1}{\gamma w_{\gamma}} )\! L_{g,1}^{2}\gamma ^2( \!\frac{1}{m}\!\| x^k\!-\!1_m\!\otimes\! \bar{x}^k \| ^2\!+\!\frac{1}{m}\!\| z^k\!-\!1_m\!\otimes\! \bar{z}^k \| ^2 )
\nonumber \nonumber \\
\leqslant & ( 1-\frac{3}{2}\frac{\mu _gL_{g,1}}{\mu _g+L_{g,1}}\gamma ) \left\| \bar{z}^k-y^*\left( \bar{x}^k \right) \right\| ^2
\nonumber \\
&\!\!+\frac{2}{w_{\gamma}}\!L_{g,1}^{2}\gamma ( \frac{1}{m}\!\| x^k-1_m\!\otimes\! \bar{x}^k \| ^2+\frac{1}{m}\!\| z^k-1_m\!\otimes\! \bar{z}^k \| ^2 ),  \nonumber
\end{align}
where the first inequality follows by adding and subtracting the term of $\nabla _yG( 1_m\otimes \bar{x}^k,1_m \otimes \bar{z}^k )$ and then employing the Young's inequality and Lipschitz continuity of $\nabla _y g_i$, while the last inequality holds due to
\begin{equation}
\begin{aligned}\label{EQ-barz-G-y}
& {\left\| \bar{z}^k-\gamma J_r\nabla _yG\left( 1_m\otimes \bar{x}^k,1_m\otimes \bar{z}^k \right) -y^*\left( \bar{x}^k \right) \right\| ^2}
\\
\leqslant &( 1-2\frac{\mu _gL_{g,1}}{\mu _g+L_{g,1}}\gamma ) \left\| \bar{z}^k-y^*\left( \bar{x}^k \right) \right\| ^2
\\
&+\gamma ( \gamma -\frac{2}{\mu _g+L_{g,1}} ) \left\| J\nabla _yG\left( 1_m\otimes \bar{x}^k,1_m\otimes \bar{z}^k \right) \right\| ^2
\\
\leqslant &( 1-2\frac{\mu _gL_{g,1}}{\mu _g+L_{g,1}}\gamma ) \left\| \bar{z}^k-y^*\left( \bar{x}^k \right) \right\| ^2.
\end{aligned}
\end{equation}
The inequality \eqref{EQ-barz-G-y} employs the strong convexity of $g(x,y)$ in $y$ and Lipschitz continuity of $\nabla _y g(x,y)$ as well as the condition $\gamma\leqslant \frac{2}{\mu _g +L_{g,1}}$. As for the term $A _2^z$, it follows from the Lipschitz continuity of $y^*(x)$ w.r.t. $x$ in Proposition \ref{PRO-PHI} that
\begin{equation}\label{EQ-A2Z}
A_{2}^{z}\leqslant L_{y^*}^{2}\left\| \bar{x}^{k+1}-\bar{x}^k \right\| ^2=L_{y^*}^{2}\alpha ^2\left\| \bar{H}_{x}^{k} \right\| ^2.
\end{equation}
Then, the desired result can be derived by substituting the inequalities \eqref{EQ-A1Z} and  \eqref{EQ-A2Z} into \eqref{EQ-zbar-ystar} and knowing the fact that $\gamma w_{\gamma} \leqslant 1$ and $
{ w_{\gamma}=\frac{\mu _gL_{g,1}}{2(\mu _g+L_{g,1})}}
$. This completes the proof. {\hfill $\blacksquare$}

\subsection{Proof of Lemma \ref{LE-inner-penalty}} \label{S-p4}
Following a similar approach as in the proof of Lemma \ref{LE-inner-errors}, the term $\left\| \bar{y}^{k+1}-y^*\left( \bar{x}^{k+1};\lambda \right) \right\| ^2$ can be bounded by the following two terms:
\begin{align}
\!\!\!\left\| \bar{y}^{k+1}\!-\!y^*( \bar{x}^{k+1};\lambda ) \right\| ^2
&\leqslant \!\!( 1\!+\!\beta w_{\beta} ) \underset{A_{1}^{y}}{\underbrace{{ \| \bar{y}^{k+1}-y^*( \bar{x}^k;\lambda ) \| ^2}}} \label{EQ-ybar-ystar} \\
&\!\!+\!( 1\!+\!\frac{1}{\beta w_{\beta}} ) \underset{A_{2}^{y}}{\underbrace{\| y^*( \bar{x}^k;\lambda ) -y^*( \bar{x}^{k+1};\lambda ) \| ^2}}. \nonumber
\end{align}
For the term $A_1^y$, it follows that
\begin{equation}\label{EQ-ybar-beta}
\begin{aligned}
A_{1}^{y}=&\left\| \bar{y}^k-\beta \bar{H}_{y}^{k}-y^*\left( \bar{x}^k;\lambda \right) \right\| ^2
\\
\leqslant &\left( 1+\beta w_{\beta} \right) { \left\| \bar{y}^k-\beta \phi -y^*\left( \bar{x}^k;\lambda \right) \right\| ^2}
\\
&+( 1+\frac{1}{\beta w_{\beta}} ) 2\beta ^2\left( L_{f,1}^{2}+\lambda ^2L_{g,1}^{2} \right)  \\
&\times ( \frac{1}{m}\left\| x^k-1_m\otimes \bar{x}^k \right\| ^2+\frac{1}{m}\left\| y^k-1_m\otimes \bar{y}^k \right\| ^2 )
\\
\leqslant& ( 1-\frac{3}{2}\frac{\mu _{\lambda}L_{\lambda}}{\mu _{\lambda}+L_{\lambda}}\beta ) \left\| \bar{y}^k-y^*\left( \bar{x}^k;\lambda \right) \right\| ^2
\\
&+\frac{4}{w_{\beta}}\beta ( L_{f,1}^{2}+\lambda ^2L_{g,1}^{2} ) \\
&\times ( \frac{1}{m}\left\| x^k-1_m\otimes \bar{x}^k \right\| ^2+\frac{1}{m}\left\| y^k-1_m\otimes \bar{y}^k \right\| ^2 ),
\end{aligned}
\end{equation}
where the term $\phi$ is given by
$$
{ \phi \!=\!J_r\nabla _yF\left(\! 1_m\otimes \bar{x}^k,1_m\!\otimes\! \bar{y}^k \right) \!+\!\lambda J_r\nabla _yG\!\left( \!1_m\!\otimes \!\bar{x}^k,1_m\!\otimes\! \bar{y}^k \right) }.
$$
Additionally, the last step in \eqref{EQ-ybar-beta}   is derived by
\begin{equation}
\begin{aligned}
&{\left\| \bar{y}^k-\beta \phi -y^*\left( \bar{x}^k;\lambda \right) \right\| ^2}
\\
\leqslant &( 1-2\frac{\mu _{\lambda}L_{\lambda}}{\mu _{\lambda}+L_{\lambda}}\beta ) \left\| \bar{y}^k-y^*\left( \bar{x}^k;\lambda \right) \right\| ^2
\\
&-( \frac{2}{\mu _{\lambda}+L_{\lambda}}-\beta ) \beta \left\| \phi \right\| ^2
\\
\leqslant& ( 1-2\frac{\mu _{\lambda}L_{\lambda}}{\mu _{\lambda}+L_{\lambda}}\beta ) \left\| \bar{y}^k-y^*\left( \bar{x}^k;\lambda \right) \right\| ^2,
\end{aligned}
\end{equation}
where we utilize the strong convexity and gradient Lipschitz continuity of the term $f_i(x,y)+\lambda g_i(x,y)$  in $y$, as well as the condition $
\beta \leqslant\min \left\{ {\frac{2}{\mu _{\lambda}+L_{\lambda}},\frac{\mu _{\lambda}+L_{\lambda}}{2\mu _{\lambda}L_{\lambda}}} \right\}
$, which   is imposed by the  constraint $\beta\leqslant\min \left\{ \frac{1}{\lambda L_{g,1}},{ \frac{1}{2(L_{f,1}+\lambda L_{g,1})}+\frac{1}{\lambda \mu_g} } \right\}$. For the term $A_2^y$, it follows from the Lipsthiz continuity of $y^*(x;\lambda)$ in Proposition \ref{PRO-penalty}  that
\begin{equation}\label{EQ-ystar-lambda}
A_{2}^{z}\leqslant L_{y^*,\lambda}^{2}\left\| \bar{x}^{k+1}-\bar{x}^k \right\| ^2=L_{y^*,\lambda}^{2}\alpha ^2\left\| \bar{H}_{x}^{k} \right\| ^2.
\end{equation}
Then, combining the inequalities \eqref{EQ-ybar-ystar}-\eqref{EQ-ystar-lambda}, the desired result can be derived.
This completes the proof. {\hfill $\blacksquare$}

\subsection{Proof of Lemma \ref{LE-consensus errors}} \label{S-p5}
The proof is based on  the average system of the recursions in \eqref{EQ-alg-a}-\eqref{EQ-alg-c}, which are split into the following three components:

\textbf{i)}  By the recursion \eqref{EQ-alg-a} and its average system $\bar{x}^{k+1}=\bar{x}^{k}-\alpha \bar{H}_x^k$, we can derive the equation  $
x^{k+1}-1_m\otimes \bar{x}^{k+1}=\left( \mathcal{W}_n -\mathcal{J}_n \right) \left( x^k-1_m\otimes \bar{x}^k \right) -\alpha \left( H_{x}^{k}-1_m\otimes \bar{H}_{x}^{k} \right)
$ with $
\mathcal{J}_n =\frac{1_m1_{m}^{T}}{m}\otimes I_n
$. Then, by taking the squared norm on both sides of the above equation, the consensus errors  $\left\| x^{k+1}-1_m\otimes \bar{x}^{k+1} \right\| ^2$ can be bounded by:
\begin{align}
&\left\| x^{k+1}-1_m\otimes \bar{x}^{k+1} \right\| ^2 \label{EQ-x-bar-k}
\\
\leqslant& ( 1\!+\!u ) \rho\left\| x^k-1_m\otimes \bar{x}^k \right\| ^2\!+\!(1+\!\frac{1}{u}) \alpha ^2\left\| H_{x}^{k}-1_m\otimes \bar{H}_{x}^{k} \right\| ^2 \nonumber
\\
\leqslant& ( 1\!-\!\frac{1-\rho}{2} ) \!\left\| x^k\!-\!1_m\otimes \bar{x}^k \right\| ^2 \!+\!\frac{2\alpha ^2}{1-\rho}\!\left\| H_{x}^{k}-1_m\otimes \bar{H}_{x}^{k} \right\| ^2,  \nonumber
\end{align}
where the first step uses the Young's inequality with $u>0$ and the fact $\rho =\| \mathcal{W}_n -\mathcal{J}_n \| ^2=\| W - \frac{1_m 1_m^{\rm{T}}}{m} \|^2\in [0,1)$ induced by Assumption \ref{ASS-network}, and the second step follows the fact that $u$ is taken as  $
u=\frac{1-\rho}{2\rho}
$. For ease of subsequent analysis of the final term in \eqref{EQ-y-bar-k}, we introduce the following shorthand notation:  $
H_{x}^{k}( \bar{x}^k,\bar{y}^k,\bar{z}^k ) =\nabla _xF( 1_m\otimes \bar{x}^k,1_m\otimes \bar{y}^k ) +\lambda ( \nabla _xG( 1_m\otimes \bar{x}^k,1_m\otimes \bar{y}^k ) -\nabla _xG( 1_m\otimes \bar{x}^k,1_m\otimes \bar{z}^k ) )
$.  In what follows,  for the final term in \eqref{EQ-x-bar-k}, we have:
\begin{align}
&\!\left\| H_{x}^{k}\!-\!1_m\otimes \bar{H}_{x}^{k} \right\| ^2 \label{EQ-hx-hxbar}
\\
\leqslant  &3\| H_{x}^{k}\!-\!H_{x}^{k}( \bar{x}^k,y^*( \bar{x}^k ),y^*( \bar{x}^k ) ) \| ^2
\nonumber\\
\!\!&\!\!\!\!\!+\!3\|\! H_{x}^{k}( \bar{x}^k,y^*(\bar{x}^k ),y^*( \bar{x}^k ) ) \!-\!1_m\!\otimes\! [ J_nH_{x}^{k}( \bar{x}^k,y^*( \bar{x}^k )\!,y^*( \bar{x}^k ) ) ] \| ^2
\nonumber\\
&\!\!\!\!\!+\!3\| 1_m\otimes [ J_nH_{x}^{k}( \bar{x}^k,y^*( \bar{x}^k ) ,y^*( \bar{x}^k ) ) ] -1_m\otimes \bar{H}_{x}^{k} \| ^2
\nonumber\\
\leqslant &18U_{\lambda}^{2}\| x^k-1_m\otimes \bar{x}^k \| ^2+18U_{\lambda}^{2}\| y^k-1_m\otimes y^*( \bar{x}^k ) \| ^2  \nonumber \\
&+18\lambda ^2L_{g,1}^{2}\| z^k-1_m\otimes y^*( \bar{x}^k ) \| ^2\,\,+3b_{f}^{2}
\nonumber\\
\leqslant & 54U_{\lambda}^{2}\| x^k-1_m\otimes \bar{x}^k \| ^2+54U_{\lambda}^{2}\| y^k-1_m\otimes \bar{y}^k \| ^2 \nonumber \\
&\!\!+\!54\lambda ^2L_{g,1}^{2}\| z^k\!-\!1_m\!\otimes  \bar{z}^k  \| ^2 +54U_{\lambda}^{2}m\| \bar{y}^k\!-\!y^*( \bar{x}^k;\lambda ) \| ^2\!
\nonumber \\
&\!\!+54\lambda ^2L_{g,1}^{2}m\| \bar{z}^k\!-\!y^*( \bar{x}^k ) \| ^2
+\!\frac{54mC_{\rm in}^{2}U_{\lambda}^{2}}{\lambda ^2}\,\,+3mb_{f}^{2}, \nonumber
\end{align}
where  the first inequality follows from the Cauchy–Schwarz inequality; the second inequality is obtained by applying the Lipschitz continuity of $\nabla_x f_i$ and $\nabla_x g_i$ to the first two terms, and using Assumption \ref{ASS-heterogeneity} to bound the last term; the last inequality uses the relation that $
\| y^k-1_m\otimes y^*( \bar{x}^k ) \| ^2\leqslant 3\| y^k-1_m\otimes \bar{y}^k \| ^2+3m\| \bar{y}^k-y^*( \bar{x}^k;\lambda ) \| ^2+3m\| y^*( \bar{x}^k;\lambda ) -y^*( \bar{x}^k ) \| ^2
$ and the result in Lemma \ref{LE-inner-penalty-error} that $
{\| y^*( x ) -y^*( x;\lambda ) \|^2 \leqslant \frac{C_{\rm in}^2}{\lambda^2}}
$. Then, combining \eqref{EQ-x-bar-k} and \eqref{EQ-hx-hxbar} yields \eqref{EQ-consensus-x}.

\textbf{ii)} Similarly to the derivation in \eqref{EQ-x-bar-k}, combining the recursion \eqref{EQ-alg-b} and its average dynamic, we have the following evolution for the consensus errors of   $\{y_i^k\}$:
\begin{align}
&\| y^{k+1}-1_m\otimes \bar{y}^{k+1} \| ^2 \label{EQ-y-bar-k}
\\
\leqslant&( 1-\frac{1-\rho}{2} ) \| y^k \! - \!1_m\otimes \bar{y}^k \| ^2\!+\!\frac{2}{1-\rho}\beta ^2\left\| H_{y}^{k} \!- \!1_m\otimes \bar{H}_{y}^{k} \right\| ^2.  \nonumber
\end{align}
In what follows, we analyze the last term of \eqref{EQ-y-bar-k}.
Letting
$
H_{y}^{k}( \bar{x}^k,\bar{y}^k ) \triangleq \nabla _yF( 1_m\otimes \bar{x}^k,1_m\otimes \bar{y}^k ) +\lambda \nabla _yG( 1_m\otimes \bar{x}^k,1_m\otimes \bar{y}^k )
$ and employing the Cauchy–Schwarz
inequality, the last term of \eqref{EQ-y-bar-k} can be bounded by:
\begin{align}
&\left\| H_{y}^{k}-1_m\otimes \bar{H}_{y}^{k} \right\| ^2 \label{EQ-bar-hy}
\\
\leqslant & 3\left\| H_{y}^{k}-H_{y}^{k}\left( \bar{x}^k,y^*\left( \bar{x}^k \right) \right) \right\| ^2
\nonumber \\
&+3\left\| H_{y}^{k}\left( \bar{x}^k,y^*\left( \bar{x}^k \right) \right) -1_m\otimes \left( J_rH_{y}^{k}\left( \bar{x}^k,y^*\left( \bar{x}^k \right) \right) \right) \right\| ^2
\nonumber \\
&+3\| 1_m\otimes ( J_rH_{y}^{k}( \bar{x}^k,y^*( \bar{x}^k ) ) ) -1_m\otimes \bar{H}_{y}^{k} \| ^2
\nonumber \\
\leqslant & 12{U_{\lambda}^2} ( \| x^k-1_m\otimes \bar{x}^k \| ^2+\| y^k-1_m\otimes y^*( \bar{x}^k ) \| ^2 )  \nonumber \\
&+6m\left( b_{f}^{2}+\lambda ^2b_{g}^{2} \right)
\nonumber \\
\leqslant & 36{U_{\lambda}^2}( \| x^k-1_m\otimes \bar{x}^k \| ^2+\| y^k-1_m\otimes \bar{y}^k \| ^2)  \nonumber \\
&\!\!\!\!\!+ \!36{U_{\lambda}^2} m\| \bar{y}^k\!-\!y^*( \bar{x}^k;\lambda ) \| ^2 \!+ \!36 \frac{m{U_{\lambda}^2}{ C_{\rm in}^{2}}}{{\lambda ^2}}\!\!+\!6m( \!b_{f}^{2}\!+\!\lambda ^2b_{g}^{2} ), \nonumber
\end{align}
where the second inequality uses the Lipschitz continuity of $\nabla_y f_i$ and  $\nabla_y g_i$ as well as Assumption \ref{ASS-heterogeneity}; the last inequality uses the inequality $
\| y^k-1_m\otimes y^*( \bar{x}^k ) \| ^2\leqslant 3\| y^k-1_m\otimes \bar{y}^k \| ^2+3m\| \bar{y}^k-y^*( \bar{x}^k;\lambda ) \| ^2+3m\left\| y^*( \bar{x}^k;\lambda ) -y^*( \bar{x}^k ) \right\| ^2
$ and  Lemma \ref{LE-inner-penalty-error}.  Then, substituting \eqref{EQ-bar-hy} into \eqref{EQ-y-bar-k} yields \eqref{EQ-consensus-y}.

\textbf{iii)} Similarly to the derivation in \eqref{EQ-x-bar-k},  by the recursion \eqref{EQ-alg-b}, we have the following  evolution  for the consensus errors of   $\{z_i^k\}$:
\begin{align}
&\| z^{k+1}-1_m\otimes \bar{z}^{k+1} \| ^2 \label{EQ-z-bar-k}
\\
\leqslant&( 1-\frac{1-\rho}{2} ) \| z^k-1_m\otimes \bar{z}^k \| ^2+\frac{2\gamma ^2}{1-\rho}\| H_{z}^{k}-1_m\otimes \bar{H}_{z}^{k} \| ^2,  \nonumber
\end{align}
where the last term can be bounded by:
\begin{align}
&\left\| H_{z}^{k}-1_m\otimes \bar{H}_{z}^{k} \right\| ^2 \label{EQ-bar-hz}
\\
\leqslant& 3\left\| H_{z}^{k}-\nabla _yG( 1_m\otimes \bar{x}^k,1_m\otimes y^*( \bar{x}^k ) ) \right\| ^2
\nonumber\\
&+\!3\left\| \begin{array}{c}
	\nabla _yG( 1_m\otimes \bar{x}^k,1_m\otimes y^*( \bar{x}^k ) )\\
	-1_m\otimes [ J_r\nabla _yG( 1_m\otimes \bar{x}^k,1_m\otimes y^*( \bar{x}^k ) ) ]\\
\end{array} \right\| ^2
\nonumber \\
&+\!3\left\| 1_m\otimes [J_r{\nabla _yG( 1_m\!\otimes\! \bar{x}^k,1_m\!\otimes \!y^*(\! \bar{x}^k \!) )}]\!-\!1_m\!\otimes\! \bar{H}_{z}^{k} \right\| ^2
\nonumber \\
\leqslant & 12L_{g,1}^{2}( \left\| x^k-1_m\otimes \bar{x}^k \right\| ^2+\left\| z^k-1_m\otimes \bar{z}^k \right\| ^2 )
\nonumber \\
&+\!12L_{g,1}^{2}m\left\| \bar{z}^k-y^*( \bar{x}^k ) \right\| ^2+3mb_{g}^{2} \nonumber,
\end{align}
where the last inequality uses the Lipschitz continuity of $\nabla_y g_i$ and Assumption \ref{ASS-heterogeneity} and employs the relation that
$
\left\| z^k-1_m\otimes y^*\left( \bar{x}^k \right) \right\| ^2\leqslant 2\left\| z^k-1_m\otimes \bar{z}^k \right\| ^2+2m\left\| \bar{z}^k-y^*\left( \bar{x}^k \right) \right\| ^2
$. Then, substituting \eqref{EQ-bar-hz} into \eqref{EQ-z-bar-k} yields \eqref{EQ-consensus-z}.  This completes the proof. {\hfill $\blacksquare$}

\section*{References}
\bibliography{ref}
\bibliographystyle{IEEEtran}

\vspace{-0.5cm}
\begin{IEEEbiography}[{\includegraphics[width=1in,height=1.25in,clip,keepaspectratio]{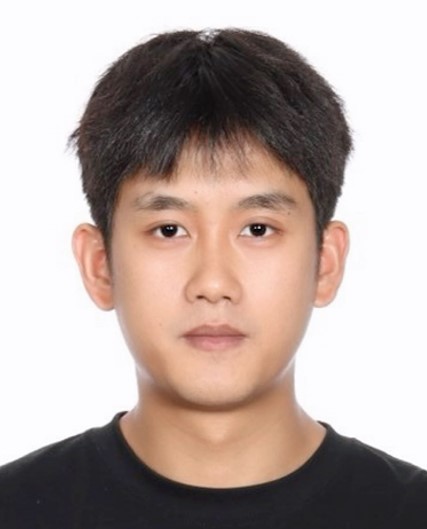}}]{Youcheng Niu} received the  B.S. degree in electronic information and engineering  and the M.E.  degree in information and communication engineering from Southwest University, Chongqing,
China, in 2020 and 2022, respectively. He is currently working toward the Ph.D degree in Zhejiang University.
His current research interests include distributed optimization and control, machine learning, and game theory.
\end{IEEEbiography}

\vspace{-0.7cm}
\begin{IEEEbiography}
[{\includegraphics[width=1in,height=1.25in,clip,keepaspectratio]{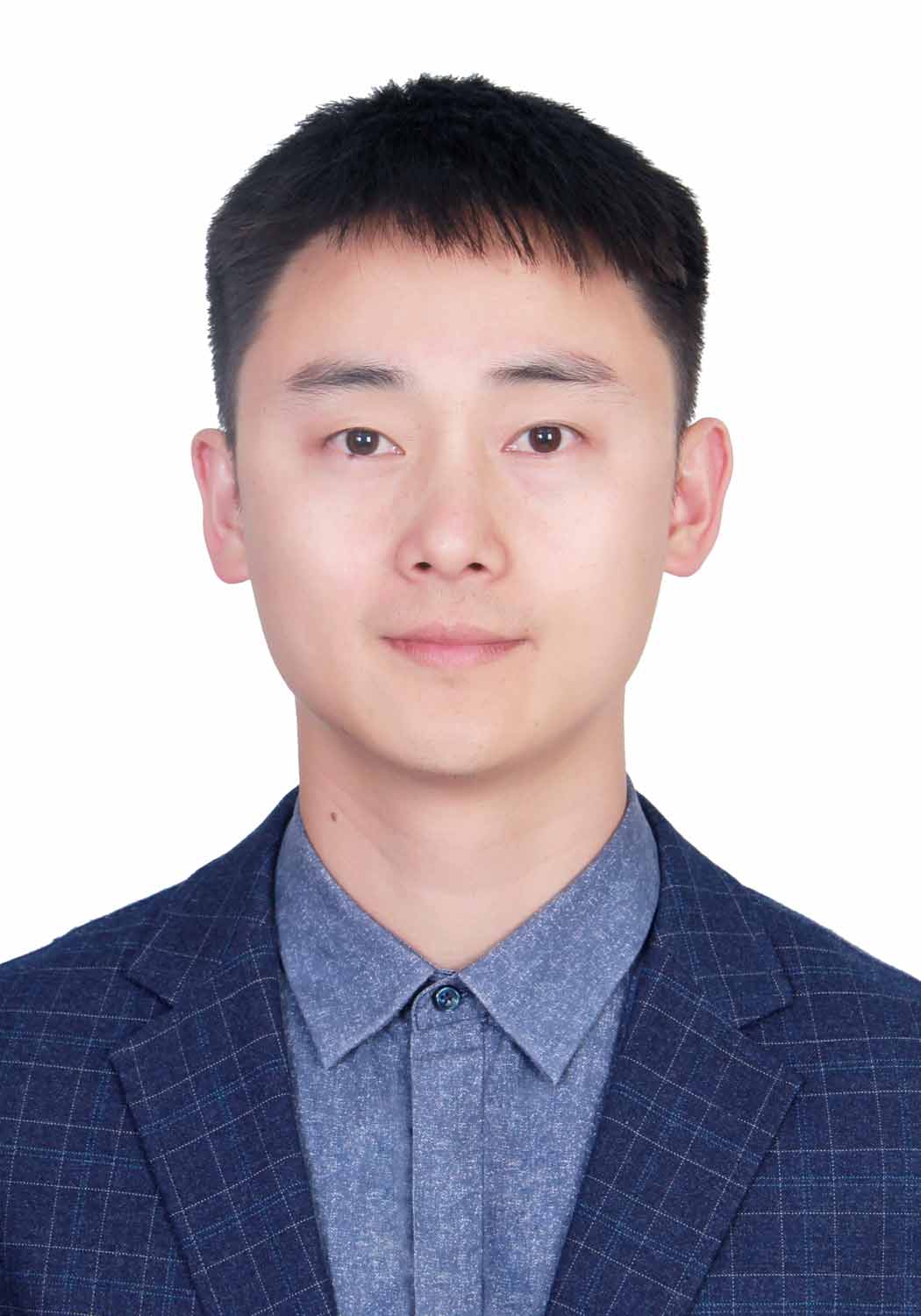}}]{Jinming Xu}
received the B.S. degree in Mechanical Engineering from Shandong University, China, in 2009 and the Ph.D. degree in Electrical and Electronic Engineering from Nanyang Technological University (NTU), Singapore, in 2016.
From 2016 to 2017, he was a Research Fellow at the EXQUITUS center, NTU; he was also a postdoctoral researcher in the Ira A. Fulton Schools of Engineering, Arizona State University, from 2017 to 2018, and the School of Industrial Engineering, Purdue University, from 2018 to 2019, respectively. In 2019, he joined Zhejiang University, China, where he is currently a Professor with the College of Control Science and Engineering. His research interests include distributed optimization and control, machine learning and network science. He has published over 50 peer-reviewed papers in prestigious journals and leading conferences. He has been the Associate Editor of IEEE Transactions on Signal and Information Processing over Networks.
\end{IEEEbiography}

\begin{IEEEbiography}
[{\includegraphics[width=1in,height=1.25in,clip,keepaspectratio]{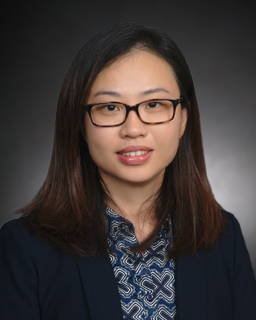}}]{Ying Sun}
received the B.E. degree in electronic information from the Huazhong University of Science and Technology, Wuhan, China,
in 2011, and the Ph.D. degree in electronic and computer engineering from
The Hong Kong University of Science and Technology in 2016. She was a
Post-Doctoral Researcher with the School of Industrial Engineering, Purdue
University, from 2016 to 2020. She is currently an Assistant Professor with the
Department of Electrical Engineering, The Pennsylvania State University. Her
research interests include statistical signal processing, optimization algorithms,
and machine learning. She was a recipient of the 2020 IEEE Signal Processing
Society Young Author Best Paper Award. She was a co-recipient of the Student
Best Paper Award from the IEEE International Workshop on Computational
Advances in Multi-Sensor Adaptive Processing (CAMSAP) 2017.
\end{IEEEbiography}

\begin{IEEEbiography}
[{\includegraphics[width=1in,height=1.25in,clip,keepaspectratio]{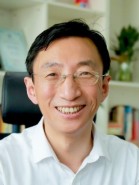}}]{Li Chai}
received the B.S. degree in applied mathematics and the M.S. degree in control science and engineering from
Zhejiang University, Hangzhou, China, in 1994
and 1997, respectively, and the Ph.D. degree in
electrical engineering from the Hong Kong University of Science and Technology, Hong Kong,
in 2002.
From 2002 to 2007, he was with Hangzhou
Dianzi University, China. He worked as a Professor with Wuhan University of Science and Technology, China, from 2008 to 2022. In 2022, he joined Zhejiang University, China, where he is currently a Full Professor with the College of Control Science and Engineering. He has been a Postdoctoral Researcher or a
Visiting Scholar with Monash University, Newcastle University, Australia,
and Harvard University, USA.

His research interests include distributed optimization, filter banks, graph signal processing, and networked control systems. He has published over 100 fully refereed papers in prestigious journals and leading conferences.
Dr. Chai is the recipient of the Distinguished Young Scholar of the
National Science Foundation of China. He is the Expert of State Council
Special Allowance, China. He has been the Associate Editor of IEEE
TRANSACTIONS ON CIRCUIT AND SYSTEMS II: EXPRESS BRIEFS, CONTROL
AND DECISION.
\end{IEEEbiography}

\begin{IEEEbiography}
[{\includegraphics[width=1in,height=1.25in,clip,keepaspectratio]{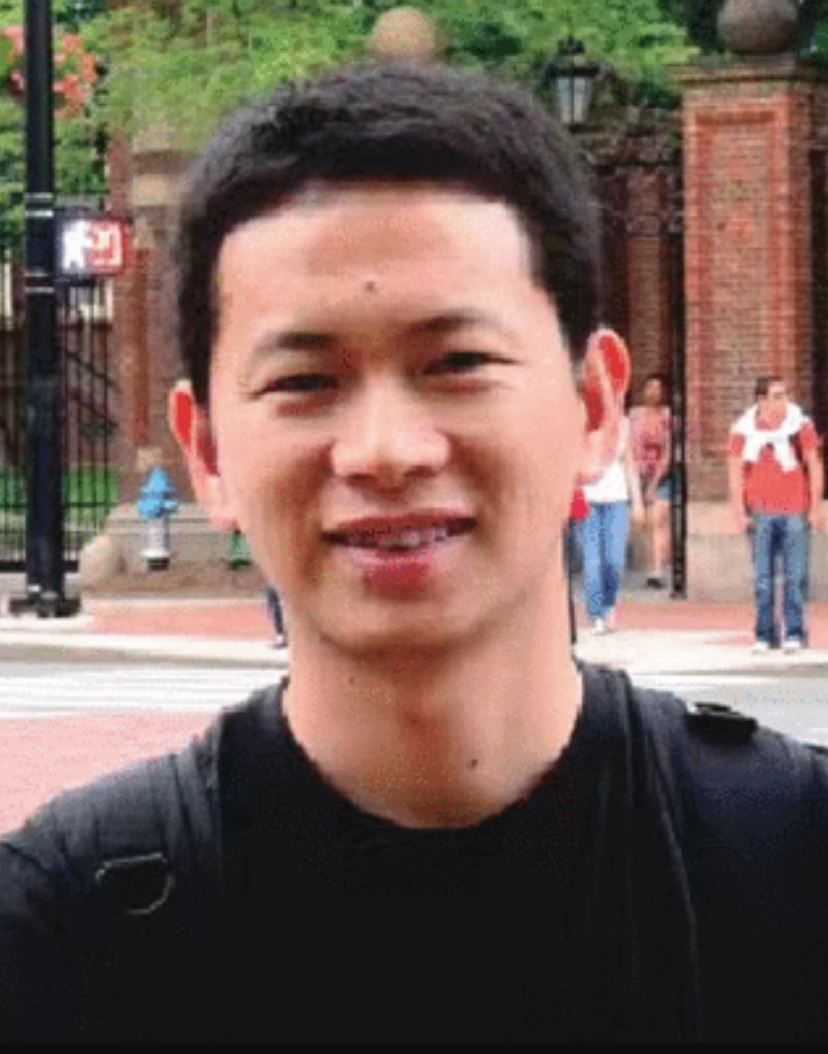}}]{Jiming Chen}(Fellow, IEEE)
 received the Ph.D.
degree in control science and engineering from
Zhejiang University, Hangzhou, China, in 2005.
He is currently a Professor with the Department of Control Science and Engineering, and the Vice Dean with the Faculty of Information Technology,
Zhejiang University. His research interests include IoT, networked control, and wireless networks. He serves on the editorial boards of multiple IEEE Transactions, and the General Co-Chairs for IEEE RTCSA’19, IEEE Datacom’19, and IEEE PST’20. He was a recipient of the 7th IEEE ComSoc Asia/Pacific Outstanding Paper Award, the JSPS Invitation Fellowship, and the IEEE ComSoc AP Outstanding
Young Researcher Award. He is an IEEE VTS Distinguished Lecturer. He is
a fellow of the CAA.
\end{IEEEbiography}



\end{document}